
\baselineskip=20pt
\magnification=\magstep1

\font\germ=eufm10
\def\ge{\hbox{\germ g}}

\def\brho{\overline{\rho}}
\def\ce#1#2{{\cal E}(#1,#2)}
\def\ow#1#2{{\cal P}(#1,#2)}
\def\owo#1#2{{\cal P}_0(#1,#2)}
\def\vp{\varphi}
\def\De{\Delta}
\def\La{\Lambda}
\def\dow{\hat{E}}
\def\rc{M\cap N'}
\def\hrc{M_1 \cap N'}
\def\ma#1#2{\sigma_{#1}^{#2}}
\def\En#1{End(#1)}
\def\Se#1{Sect(#1)}
\def\bR#1{\overline{R}_{{#1}}}
\def\cH{{\cal H}}
\def\cK{{\cal K}}
\def\cHr{\cH_\rho}
\def\rx{\rho_\xi}
\def\rbx{\rho_{\overline{\xi}}}
\def\Rx{R_\xi}
\def\bRx{\bR \xi}
\def\Hr{\cH_\rho}
\def\Hx{\cH_\xi}
\def\Hbx{\cH_{\overline{\xi}}}
\def\Kx{{\cal K}_\xi}
\def\dx{d(\xi)}
\def\Ex{E_\xi}

\def\cx{c_\xi}
\def\bcx{\overline{c}_\xi}
\def\ax{a_\xi}
\def\Ax{A_\xi}
\def\bax{\overline{a}_\xi}
\def\Abx{A_{\overline{\xi}}}
\def\bV{\overline{V}}
\def\nx{n_\xi}
\def\zx{z_\xi}
\def\zbx{z_{\overline{\xi}}}
\def\axi{a_{\xi,i}}
\def\Vxi{V(\xi)_i}
\def\xxi{x(\xi)_i}
\def\pxi{p_{\xi,i}}
\def\cK{{\cal K}}
\def\Kx{{\cal K}_\xi}
\def\Kxo{\Kx^\perp}
\def\mx{m_\xi}
\def\cC{{\cal C}}
\def\hG{\widehat{G}}
\def\Vpi{V(\pi)_i}
\def\cA{{\cal A}}
\def\dA{\widehat{\cA}}
\def\LA{L^2(\cA)}
\def\dd{\widehat{\delta}}

\def\dk{\widehat{\kappa}}
\def\dV{\widehat{V}}
\def\cB{{\cal B}}
\def\cC{{\cal C}}
\def\tB{{\widetilde{B}}}
\def\Oh{\Omega_h}
\def\dh{\hat{h}}
\def\Odh{\Omega_{\dh}}
\def\tb{\tilde{b}}

\font\ftitle=cmbx10 scaled\magstep1


\magnification1200
\centerline{\ftitle A Galois Correspondence}
\centerline{\ftitle for Compact Groups of Automorphisms}
\centerline{\ftitle of von Neumann Algebras}
\centerline{\ftitle with a Generalization to Kac Algebras}
\bigskip\bigskip
\centerline{ Masaki Izumi\footnote{$^*$}{Miller Research Fellow.}}
\centerline{Department of Mathematics, University of California, Berkeley, CA
94720, USA}
\smallskip
\centerline{ Roberto Longo\footnote{$^{**}$}{Supported in part by
MURST and CNR-GNAFA.}}
\centerline{Dipartimento di
Matematica, Universit\`a  di Roma ``Tor Vergata''}
\centerline{Via della Ricerca Scientifica, I-00133 Roma, Italy}
\smallskip
\centerline{ Sorin Popa\footnote{$^{***}$}{Supported in part by NSF Grant
DMS-9206984.}}
\centerline{Department of Mathematics, University of California,
Los Angeles, CA 90024, USA}
\bigskip
\bigskip
\bigskip
\bigskip
\bigskip

\noindent{\it Abstract.} Let $M$ be a factor with separable predual
and $G$ a compact group of automorphisms of $M$ whose action is
minimal, i.e. $M^{G^\prime}\cap M = C$, where $M^G$ denotes the
$G$-fixed point subalgebra.
Then every intemediate von Neumann algebra $M^G\subset N\subset M$
has the form $N=M^H$ for some closed subgroup $H$ of $G$. An
extension of this result to the case of actions of compact Kac
algebras on factors is also presented. No assumptions are made
 on the existence of a normal conditional expectation onto $N$.
\vfill\eject

 \noindent{\bf 1. Introduction.}

A classical theme in Operator Algebras is the Galois
correspondence between groups of automorphisms of a
von Neumann algebra and von Neumann subalgebras.

To be more specific, let $M$ be a von Neumann algebra  and to each
group $G$ of automorphisms of $M$ let associate $M^G$,
the von Neumann subalgebra of the  $G$--fixed
elements
$$
G\rightarrow M^G.\eqno (1.1)
$$
In a dual way to each von Neumann subalgebra
$N$ of $M$ we may associate the group $G_N$ of the
automorphisms of $M$ leaving $N$ pointwise fixed
$$
N\rightarrow G_N.\eqno (1.2)
$$
These two maps are in general not one another inverse,
but restricting to (closed) subgroups of a given  group $G$ and
to intermediate von Neumann subalgebras $M^G\subset N\subset M$
they may actually become one another inverse.

Such a Galois correspondence was shown to hold by Nakamura and
Takeda [NT] and  Suzuki [Su] in the case $M$ be $II_1$--factor
and $G$ a finite group whose action on $M$ is minimal,
namely $M^{G^\prime}\cap M = \bf C$.

A different Galois correspondence, between normal closed subgroups of
a compact (minimal) group $G$ and globally $G$-invariant intermediate
von Neumann algebras, was obtained by Kishimoto [K], following
methods in the analysis of the chemical
potential in Quantum Statistical Mechanics [AHKT]. Generalizations of
this result concerning dual actions of a locally compact group $G$
were dealt by Takesaki, in case of $G$  abelian, and by Nakagami in
more generality, see [NTs].

Another kind of Galois correspondence was provided by H. Choda [Ch].
It concerns in particular the crossed product of a factor by an
outer action of a discrete group and characterizes the intermediate
von Neumann subalgebras that are crossed product by a discrete
subgroup. An important assumption here is the existence of a normal
conditional expectation onto the intermediate subalgebras.

In this paper we consider any  compact
group $G$ of automorphisms of a (separable) factor $M$, whose
action is minimal, and show that any intermediate von Neumann algebra
$M^G\subset N\subset M$ is the fixed-point algebra $N=M^H$ for some
closed subgroup $H$ of $G$, namely  the general Galois
correspondence holds in the compact minimal case. Indeed as a
corollary the two maps $(1.1)$ and $(1.2)$ are one another inverse.

A particular case of our result concerning the action
of the periodic modular group with minimal spectrum on a type
$III_\lambda$ factor, $0<\lambda<1$, has been recently obtained in
[HS].

Concerning the ingredients in our proof, we mention the spectral
analysis for compact group actions,  endomorphisms and index theory
for infinite factors,  arguments based on modular
theory, injective subfactors and averaging techniques.

We enphasise that the main step in the proof of our result is
showing  the existence of a (necessarily unique) normal  conditional
expectation of $M$ onto any intermediate subfactor between $M^G$
and $M$.

Note that in this way we also obtain a Galois correspondence for
intermediate von Neumann algebras in the case of crossed products
of factors by outer actions of discrete groups (again, without the
a priori existence of normal conditional expectation).

 This poses the following question: if $M_1\subset M_2\subset
M_3$ are von Neumann algebras such that $M_1'\cap M_3 = C$, with a
normal expectation $\epsilon : M_3\to M_1$, does there exist a normal
expectation of $M_3$ onto $M_2$? In other words, does $\epsilon$
factor through $M_2$? Beside the case dealt in this paper, we know a
(positive) answer  for example if $M_1\subset M_3$ has finite index
or if $M_3$ is semifinite, but no counter-example is known to us.

At this point we briefly comment on the superselection
structure in Particle Physics, that partly motivated
our work. As is known the group of the internal symmetries in a
Quantum Field Theory is
the dual of the tensor C$^*$-category defined by the superselection
sectors [DR]. Our result classifies the  extensions of the net of
the observable algebras made up by field operators.
An analysis of further aspects of this structure goes
beyond the purpose of our paper. However we notice that in
low dimensional Quantum Field Theory the internal symmetry is
realized by a more general, not yet understood, quantum object
and this suggests to be of interest to extend our result to
a wider class of ``quantum groups''.

We take here a first step in this direction providing a version
of our result in the context of actions of compact Kac algebras on
factors that turns out to be new even in the finite-dimensional
case. This is included in our last section.

\vskip15pt
\noindent {\bf 2. Preliminaries.}

Throughout this paper, von Neumann algebras have separable preduals.

\noindent {\it 2-1. Operator valued weights and basic construction.}
For the theory of operator valued weights and basic construction,
our standard references are [H1][H2][Ko1].

Let $M \supset N$ be an inclusion of von Neumann algebras.
We denote by $\ow MN$, $\ce MN$ the set of normal semifinite faithful,
(abbreviated as n.s.f.), operator valued weights, and that of normal faithful
conditional expectations respectively.
We denote by $\owo MN$ the set of $T\in \ow MN$ whose restriction to
$M\cap N'
$ is semifinite, (such $T$ is called regular in [Y]).
Note that $\owo MN$ is either empty or $\ow MN$ [H2, Theorem 6.6].
For $T\in \ow MN$, we use the following standard notations:
$$n_T=\{ x\in M; T(x^*x)< \infty \},$$
$$m_T=n_T^* n_T.$$
For a n.f.s. weight $\vp$ on $M$, $H_\vp$ and $\La_\vp$ denote the
GNS Hilbert space and the canonical injection
$\La_\vp : n_\vp \longrightarrow H_\vp$.

For $M\supset N$ with $E\in \ce MN$, we fix a faithful normal state $\omega$
on $N$ and set $\vp:=\omega \cdot E$.
We regard $M$ as a concrete von Neumann algebra acting on $H_\vp$.
Let $e_N$ be the Jones projection defined by $e_N\La_\vp(x)=\La_\vp(E(x))$,
which does not depend on $\omega$ but only on the natural cone of
$H_\vp$ [Ko2, Appendix].
The basic extension of $M$ by $E$ is the von Neumann algebra generated by $M$
and $e_N$, which coincides with $J_M N' J_M$, where $J_M$ is the modular
conjugation for $M$.
For $x\in B(H_\vp)$, we set $j(x)=J_M x^* J_M$.
The dual operator valued weight $\dow \in \ow {M_1}M$ of $E$ is defined by
$j\cdot E^{-1} \cdot j|_{M_1}$, where $E^{-1}\in \ow {N'}{M'}$ is
characterized by
spatial derivatives [C1]:
$${d (\psi \cdot E) \over d \vp'}={d \psi \over d (\vp' \cdot E^{-1})}, \quad
\psi \in {\cal P}(N,{\bf C}), \enskip \vp' \in {\cal P}(M',{\bf
C}).$$ Since $\dow$ satisfies $\dow(e_N)=1$ [Ko1, Lemma 3.1],
$Me_N M \subset m_{\dow}$.
In [Ko1], Kosaki defined the index of $E$ by $Ind \;E=E^{-1}(1)$ in the case
where $M$ and $N$ are factors, which is known to coincide with
the probabilistic index defined in [PP1].

First, we consider a Pimsner-Popa push down lemma in our setting
(c.f. [PP1]).

\proclaim Lemma 2.1.
Let $M$ be a von Neumann algebra and $\vp$ a n.f.s. weight on $M$.
Suppose ${\cal A}$ is  a *-subalgebra of $n_\vp^* \cap n_\vp$ which is dense
in $M$ in weak topology, and globally invariant under the modular automorphism
group.
Then $\La_\vp ({\cal A})$ is dense in $H_\vp$. \par

\noindent {\it Proof.}
Let $p$ be the projection onto the closure of $\La_\vp({\cal A})$.
Then $p \in {\cal A}'=M'$.
Thanks to $\sigma_t^\vp({\cal A})={\cal A}$, $p$ commutes with
$\De_\vp^{it}$, and consequently we have
$\De_\vp^{1\over 2}p \supset p\De_\vp^{1\over 2}$.
Since ${\cal A}\subset n_\vp \cap n_\vp^*$, $\La_\vp(x)$, $x\in {\cal A}$ is
in the domains of $S_\vp$ and $\De_\vp^{1\over 2}$.
Thus we get the following:
$$J_\vp \La_\vp(x)=J_\vp S_\vp \La_\vp(x^*)=\De_\vp^{1\over 2}\La_\vp(x^*)
=p\De_\vp^{1\over 2}\La_\vp(x^*)=pJ_\vp \La_\vp(x).$$
This means that $p$ commutes with $J_\vp$, and $p\in M\cap M'$.
So we get $\La_\vp((1-p)x)=0$ for $x\in {\cal A}$.
Since $\vp$ is faithful, this implies $(1-p)x=0$, which shows $p=1$ because
${\cal A}$ is dense in $M$ in weak topology. Q.E.D.

\proclaim Proposition 2.2 (Push down lemma).
Let $M \supset N$ be an inclusion of factors with $E \in \ce MN$,
and $M_1$ be the basic extension of $M$ by $E$.
Then for all $x \in n_{\dow}$, $e_N\dow(e_Nx)=e_Nx$ holds. \par

\noindent {\it Proof.} Let $\vp$ be as above and $\vp_1 =\vp \cdot \dow$.
Then $e_Nx$, and $e_N\dow(e_Nx)$ belong to $n_{\vp}$.
So we get the following:
$$\eqalign{
&||\La_{\vp_1}(e_Nx)-\La_{\vp_1}(e_N\dow(e_Nx))||^2 \cr
&=\vp_1(x^*e_Nx)-\vp_1(x^*e_N\dow(e_Nx))-\vp_1(\dow(e_Nx)^*e_Nx)
+\vp_1(\dow(e_Nx)^*e_N\dow(e_Nx)) \cr
&=||\La_{\vp_1}(e_Nx)||^2-||\La_{\vp_1}(e_N\dow(e_Nx))||^2.}$$
So, we can define a bounded operator $V$ on $e_NH_{\vp_1}$ by
$$Ve_N \La_{\vp_1}(x)=\La_{\vp_1}(e_N\dow(e_Nx)), \quad x\in n_{\dow}. $$
By simple computation, one can show that $V$ is identity on
$e_N\La_{\vp_1}(Me_NM)$.
So to prove the statement, it suffices to show that $\La_{\vp_1}(Me_NM)$ is
dense in $H_{\vp_1}$.
We set ${\cal A}=Me_NM$ and show that ${\cal A}$ satisfies the assumption of
the previous lemma.
Indeed, since $M_1$ is the weak closure of $Me_NM+M$, the weak closure of
${\cal A}$ is a closed two-sided ideal of $M_1$, and coincides
with $M_1$.
From the definition of $\vp_1$, we have $\sigma_t^{\vp_1}(Me_NM)=
\sigma_t^\vp(M)\sigma_t^{\vp_1}(e_N)\sigma_t^\vp(M)=M\sigma_t^\vp(e_N)M$.
Thanks to $j\cdot E^{-1} \cdot j=(j\cdot E \cdot j)^{-1}$ [Ko1, Lemma 1.3],
we get
$${d\vp_1 \over d (\omega \cdot j)}
={d(\vp\cdot (j\cdot E\cdot j)^{-1}) \over d (\omega \cdot j)}
={d \vp \over d (\omega \cdot E \cdot j)}={d\vp \over d(\vp \cdot j)}
=\De_\vp.$$
Since $\De_\vp$ commutes with $e_N$, we get $\sigma_t^{\vp_1}(e_N)=e_N$.
Q.E.D.

\noindent {\it Remark 2.3.} In general, $e_NM_1$ is strictly larger than
$e_N n_{\dow}$.
Indeed, suppose that $M$, $N$ are type III factors and $e_NM_1=e_Nn_{\dow}$.
Then there is an isometry $v\in n_{\dow}$ with $vv^*=e_N$, that implies
$1=v^*v\in m_{\dow}$ and $Ind \;E<\infty$.
This also means that $e_NM$ is not necessarily closed in weak topology because
of $e_NM_1=\overline{e_NM}^w$.

The following is a generalization of the abstact characterization of
the basic extension in [PP2] to the
infinite index case (see also [HK]).

\proclaim Lemma 2.4.
Under the same assumption, assume that $R$ is a factor including $M$ and
satisfying the following:
\item {(i)} There is a projection $e\in R$ such that $R$ is generated by
$e$ and $M$, and $exe=E(x)e$ holds for $x\in M$.
\item {(ii)} There is $T\in \ow RM$ satisfying $T(e)=1$, and
$e\in (R\cap N')_{E\cdot T}$. \par
\noindent {\it
Then there is an isomorphism $\pi :M_1 \longrightarrow R$ satisfying
$\pi|_M=id_M$, $\pi(e_N)=e$, and $T\cdot \pi =\pi \cdot \dow$.}

\noindent {\it Proof.} Let $\psi=\vp \cdot T$.
For the same reason as before, $\La_\psi (MeM)$ is dense in $H_\psi$.
So we can define a surjective isometry $U:H_{\vp_1} \longrightarrow H_\psi$
and an isomorphism $\pi :M_1 \longrightarrow R$ by
$$U\La_{\vp_1}(\sum x_ie_Ny_i)=\La_{\psi}(\sum x_iey_i),
\quad x_i, y_i, \in M,$$
$$\pi(x)=UxU^*, \quad x\in M_1.$$
Clearly, $\pi$ satisfies $\pi|_M=id_M$, $\pi(e_N)=e$.
Thanks to $\sigma_t^{E\cdot T}(e)=e$, the modular automorphism groups of
$\vp_1$ and $\psi \cdot \pi$ coincide on $Me_NM$ (and on $M_1$).
Since $M_1$ is a factor, this implies that $\vp_1$ is a scalar multiple of
$\psi \cdot \pi$, and consequently that $\dow$ is a scalar multiple of
$\pi^{-1} \cdot T \cdot \pi$ [H2, Lemma 4.8].
From $\dow(e_N)=1$ and $T(e)=1$, we get the result. Q.E.D.

The following may be a folklore for specialists.
However, since the authors can not find it in the literature, we give a proof.

\proclaim Lemma 2.5.
Let $M\supset N$ be an inclusion of von Neumann algebras (not necessarily
with separable pre-dual).
Then, there is a unique central projection $z$ of $\rc$ satisfying the
following two conditions:
\item {(i)} $\owo {pMp}{pN}=\emptyset$ holds for every projection
$p\in \rc$, $p\leq 1-z$.
\item {(ii)} $\owo {zMz}{zN}=\ow {zMz}{zN}$. \par
\noindent {\it
Moreover, if $\ow MN$ is not empty, then $(1-z)(M \cap N') \cap m_T=\{0\}$,
$z\in (\rc)_T$, and $T|_{z(\rc)}$ is semifinite for every $T\in \ow MN$}
\par

To prove the lemma, we need the following.

\proclaim Lemma 2.6. The following hold.
\item {(i)} Let $\{p_i\}_{i\in I}\subset \rc$ be a family of mutually
orthogonal projections, and $p=\sum p_i$.
If $\owo {p_iMp_i}{p_iN}\neq \emptyset$ for every $i\in I$,
then $\owo {pMp}{pN}\neq \emptyset$.
\item {(ii)} Let $p\in \rc$ be a projection.
If $\owo MN\neq \emptyset$, then $\owo {pMp}{pN}\neq \emptyset$.
\item {(iii)} Let $p \in \rc$ be a projection satisfying
$\owo {pMp}{pN}\neq \emptyset$, and $c(p)$ the central support of $p$ in
$\rc$.
Then $\owo {c(p)Mc(p)}{c(p)N}\neq \emptyset$,
\item {(iv)} Let $\{p_i\}_{i \in I}\subset \rc$ be a family of
projections, and $p_0=\vee p_i$.
If $\owo {p_iMp_i}{p_iM}\neq \emptyset$ for every $i\in I$,
then $\owo {p_0Mp_0}{p_0N}\neq \emptyset$. \par

\noindent{\it Proof.}
(i): This follows from the following easy facts:
$$\owo {pMp}{\oplus p_iMp_i}\neq \emptyset,\quad
\owo {\oplus p_iMp_i}{\oplus p_iN}\neq \emptyset,\quad
\owo {\oplus p_iN}{pN}\neq \emptyset.$$
(ii): Since $\owo MN \neq \emptyset$, there is a separating family of normal
conditional expectations from $M$ to $N$ $\{E_\alpha\}$ [Ha2, Theorem 6.6].
Then $\{E_\alpha (p\cdot p)p\}$ is a separating family of bounded normal
operator valued weights from $pMp$ to $pN$, and
$\owo {pMp}{pN}\neq \emptyset$.
(iii): Let $\{e_j\}\in \rc$ be a family of projections satisfying
$p\succ e_j$, $\sum e_j=c(p)$.
Then $\owo {e_jMe_j}{e_j N}\neq \emptyset$.
So using (i),  we get $\owo {c(p)Mc(p)}{c(p)N}\neq \emptyset$.
(iv): Let $z_i=c(p_i)$ and $z_0=\vee z_i$.
Then thanks to (i)(ii)(iii), $\owo {z_0Mz_0}{z_0N}\neq \emptyset$.
Since $z_0$ is the central support of $p_0$,
we get the statement by using (i). Q.E.D.

{\noindent}{\it Proof of Lemma 2.5.}
Let $z$ be the supremum of the projections $p\in \rc$ satisfying
$\owo {pMp}{pN}\neq \emptyset$.
Then thanks to Lemma 2.6 (iii)(iv), $z$ is a central projection satisfying
(i)(ii).
It is easy to show the uniqueness of such a projection.
If $T\in \ow MN$, $\sigma_t^T(z)$ also satisfies (i)(ii) and we get
$z\in (\rc)_T$.
This implies that $zT(z \cdot z)$ belongs to $\ow {zMz}{zN}$.
So due to [Ha2, Theorem 6.6], $T|_{z(\rc)}$ is semifinite.
Suppose $x$ is a non-zero positive element in $m_T \cap (1-z)(\rc)$.
Then there is a non-zero spectral projection $p$ of $x$ satisfying
$T(p)<\infty$.
This implies $\ce {pMp}{pN}\neq \emptyset$, that contradicts Lemma 2.6 (iii).
Q.E.D.

To analyze local structure of the inclusions obtained by basic construction in
infinite index case, we need the following:

\proclaim Lemma 2.7 ([Ko1, Proposition 4.2][Y, Corollary 28]).
Let $M \supset N$ be an inclusion of factors.
Then the following hold:
\item {(i)} Let $T \in \ow MN$, and $p\in m_T \cap (\rc)_T$ a non-zero
projection.
Then $Ind\; T_p=T(p)T^{-1}(p)$, where $T_p \in \ce {pMp}{pN}$ is defined by
$T_p(x)=pT(x)/ T(p)$, $x\in pMp$.
\item {(ii)} If $\owo MN\neq \emptyset$, $\owo {N'}{M'}\neq \emptyset$,
then $\rc$ is a direct sum of type I factors and $pMp \supset pN$ has
finite index for every finite rank projection in $\rc$. \par

\proclaim Proposition 2.8. Let $M\supset N$ be an inclusion of factors with
$E\in \ce MN$, and $M_1$ the basic extension.
Then $\hrc$ is direct sum of four subalgebras,
$$\hrc =A\oplus B_1 \oplus B_2 \oplus C,$$
satisfying the following:
\item {(i)} Each of the four subalgebras is globally invariant under
$\{\ma t{E\cdot \dow} \}$.
\item {(ii)} $j(A)=A$, $j(B_1)=B_2$, $j(B_2)=B_1$, $j(C)=C$.
\item {(iii)} $\dow |_{A\oplus B_1}$ is semifinite.
\item {(iv)} $m_{\dow} \cap (B_2 \oplus C)=\{0\}$.
\item {(v)} $A$ is direct sum of type I factors and $pM_1p \supset pN$ has
finite index for every finite rank projection $p \in A$. \par

\noindent {\it Proof.}
First, we show $j\cdot \ma t{E\cdot \dow} \cdot j=\ma t{E\cdot \dow}$ on
$\hrc$.
Indeed, for $x \in \hrc$ we get the following as in the proof of Lemma 2.3:
$$\eqalign{
j\cdot \ma t{E\cdot \dow}(j(x))
&=J_M({d(\vp \cdot \dow)\over d(\omega \cdot j)})^{it}J_M x
J_M({d(\vp \cdot \dow)\over d(\omega \cdot j)})^{-it}J_M \cr
&=J_M\De_\vp^{it}J_M x J_M\De_\vp^{-it}J_M =\De_\vp^{it} x \De_\vp^{-it}
=({d(\vp \cdot \dow)\over d(\omega \cdot j)})^{it} x
({d(\vp \cdot \dow)\over d(\omega \cdot j)})^{-it} \cr
&=\ma t{E \cdot \dow}(x).}$$
Now, let $z$ be the central projection of $\hrc$ determined by Lemma 2.5
for $M_1 \supset N$.
We set
$$A=zj(z)(\hrc), \quad C=(1-z)j(1-z)(\hrc),$$
$$B_1=zj(1-z)(\hrc), \quad   B_2=(1-z)j(z)(\hrc). $$
Then by construction (ii)(iii)(iv) hold.
Since $j$ commutes with $\ma t{E\cdot \dow}$, $j(z) \in (\hrc)_{E\cdot \dow}$,
and we get (i).
Note that for a projection $p\in \hrc$, $J_M(pM_1p)'J_M=j(p)N$,
$J_M(pN)'J_M=j(p)M_1j(p)$.
So $(pN)' \supset (pM_1p)'$ is anti-conjugate to $j(p)M_1j(p)\supset j(p)N$.
Thus thanks to Lemma 2.7, we get (v). Q.E.D.

\vskip10pt
\noindent{\it 2-2. Sectors and simple injective subfactors.}
Our basic references for the theory of sectors are [L1][L2][I1].

Let $M$ be an infinite factor.
We denote by $\En M$ and $\Se M$ the set of unital endomorphisms of $M$ and
that of sectors, which is the quotient of $\En M$ by the unitary equivalence.
Note that every element in $\En M$ is automatically normal for $M$ with
separable pre-dual.
For $\rho_1, \rho_2 \in \En M$, $(\rho_1, \rho_2)$ denotes the set of
intertwiners between $\rho_1$ and $\rho_2$, i.e.
$$(\rho_1, \rho_2)=\{v\in M;v\rho_1(x)=\rho_2(x)v, \enskip x\in M\}.$$
If $\rho_1$ is irreducible, i.e. $M\cap \rho_1(M)'={\bf C}$,
$(\rho_1, \rho_2)$ is a Hilbert space with the following inner product:
$$<V|W>1=W^*V, \quad V,W\in (\rho_1, \rho_2).$$
We define the dimension $d(\rho)$ of $\rho$ by $d(\rho)=[M:\rho(M)]_0^{1/2}$,
where $[M:\rho(M)]_0$ is the minimum index of $M \supset \rho(M)$ [Hi].
For $\rho$ with $d(\rho)<\infty$, we denote by $E_\rho$ and $\phi_\rho$
the minimal conditional expectation onto $\rho(M)$ and the standard left
inverse of $\rho$, i.e. $\phi_\rho=\rho^{-1}\cdot E_\rho$.

There are three natural operations in $\Se M$: the sum, the product and
the conjugation.
For simplicity, we denote by $\brho$ one of the representatives of
the conjugate sector
$\overline{[\rho]}$ of $[\rho]$.
When $d(\rho)$ is finite, it is known that there are two isometries
$R_\rho \in (id, \brho\rho)$, $\bR \rho \in (id, \rho\brho)$ satisfying
$$\bR \rho^*\rho(R_\rho)=R_\rho^*\brho(\bR \rho)={1\over d(\rho)}.
\eqno (2.1)$$
Although such a pair is not unique, we fix it once and forever in this paper.
Unless $\rho$ is a pseudo-real sector [L1], we can take $\bR \rho$ equal to
$R_{\brho}$.
If it is, we set $\bR \rho=-R_\rho$.

Let $_MX_M$ be a $M-M$ bimodule, and $\rho\in \En M$.
Then we define a new bimodule $_M(X_\rho)_M$ (respectively $_M(_\rho X)_M$) by
$$x\cdot \xi' \cdot y:=x\cdot \xi \cdot \rho(y)\enskip
({\rm respectively}\; x\cdot \xi' \cdot y:=\rho(x)\cdot \xi \cdot y)
\enskip x,y\in M,$$
where $\xi'=\xi$ as an element of Hilbert space $X$.
It is known that there is one-to-one correspondence between $\Se M$ and the
set of equivalence classes of $M-M$ bimodules.
The correspondence is given by $[\rho]\rightarrow [_M(L^2(M)_\rho)_M]$,
that preserves the three operations.
The conjugate sector of $[\rho]$ is characterized by
$$_M(L^2(M)_\rho)_M \simeq _M(_{\brho} L^2(M))_M.$$

Let $\phi$ be a unital normal completely positive map from $M$ to $M$.
Following Connes [C2], there is a natural way to associate a $M-M$ bimodule
with $\phi$.
Let $\Omega$ be a separating and cyclic vector of $M$.
We introduce a positive semi-definite sesquilinear form on the algebraic
tensor product $M\otimes_{alg}M$ as follows:
$$<\sum_i x_i\otimes y_i, \sum_j z_j\otimes w_j>
=\sum_{i,j}<\phi(z_j^*x_i)J_Mw_jy_i^*J_M\Omega|\Omega>.$$
We denote by $H_\phi$ the Hilbert space completion of the quotient of
$M\otimes_{alg}M$ by the kernel of the sesquilinear form, and by
$\La_\phi$ the natural map $\La_\phi :M\otimes_{alg}M \longrightarrow H_\phi$.
$H_\phi$ is naturally a $M-M$ bimodule by the following action:
$$x\cdot \La_\phi(\sum_i z_i\otimes w_i)\cdot y
=\La_\phi(\sum_ixz_i\otimes w_iy). $$
Thanks to the one-to-one correspondence stated above, there is
an endomorphism $\rho_\phi$ satisfying
$_M(H_\phi)_M\simeq  _M(_{\rho_\phi}L^2(M))_M$.
Actually, $\rho_\phi$ is Steinspring type dilation of $\phi$.
Indeed, let $W:H_\phi \longrightarrow L^2(M)$ be the intertwining surjective
isometry, and set $\xi_0=W\La_\phi(1\otimes 1)$.
Then we get
$$<\phi(x)\cdot \Omega\cdot y|\Omega>
=<x\cdot \La_\phi(1\otimes 1)\cdot y, \La_\phi(1\otimes 1)>
=<\rho_\phi(x)\cdot \xi_0\cdot y|\xi_0>.$$
We define an isometry $v$ by $v(\Omega \cdot y)=\xi_0\cdot y$.
Then by definition, $v$ commutes with the right action of $M$.
So $v$ belongs to $M$ and satisfies $\phi(x)=v^*\rho_\phi(x)v$, $x\in M$.
Note that the support of $vv^*$ in $M\cap \rho_\phi(M)'$ is $1$.
Indeed, suppose $z\in M\cap \rho_\phi(M)'$ satisfying
$z\xi_0\cdot y=0$, for all $y\in M$.
Then $z\rho_\phi(x)\xi_0\cdot y=0$ for all $x,y\in M$.
Since $\overline{\rho_\phi(M)\cdot \xi_0 \cdot M}=WH_\phi=L^2(M)$,
we get $z=0$.

Although the following statements might be found in the literature,
we give proofs for readers' convenience.

\proclaim Proposition 2.9. Let $M$ and $\phi$ be as above.
Then the following hold:
\item {(i)} Let $\sigma \in \En M$, and $v_1 \in M$ be an isometry satisfying
$\phi(x)=v_1^*\sigma(x)v_1$.
If the support of $v_1v_1^*$ in $M\cap \sigma(M)'$ is 1, then
$[\rho_\phi]=[\sigma]$.
\item {(ii)} The equivalence class of $H_\phi$ does not depend on the choice
of the cyclic separating vector $\Omega$.
\item {(iii)} Let $\mu$ be another unital normal completely positive map from
$M$ to $M$.
If there is a positive constant $c$ such that $c\mu-\phi$ is completely
positive, then $[\rho_\mu]$ contains $[\rho_\phi]$. \par

\noindent{\it Proof.} (i): Let $\xi_0$ be as before.
Then by assumption, we get the following:
$$<\sigma(x)v_1\Omega \cdot y|v_1\Omega>=<\rho_\phi(x)\xi_0\cdot y|\xi_0>,$$
$$\overline{\sigma(M)v_1\Omega\cdot M}=L^2(M).$$
So we can define a unitary $u\in M$ by
$u\sigma(x)v_1\Omega \cdot y=\rho_\phi(x)\xi_0\cdot y$, and get
$\rho_\phi(x)=u\sigma(x)u^*$.
(ii) follows from (i).
(iii): Since $c\mu-\phi$ is completely positive, we can define a bounded map
$T:H_\mu \longrightarrow H_\phi$ by
$$T\La_\mu(\sum_i x_i\otimes y_i)=\La_\phi (\sum_i x_i\otimes y_i).$$
Then $T$ is an $M-M$ bimodule map whose image is dense in $H_\phi$.
Let $T=U|T|$ be the polar decomposition of $T$.
Then $U$ is a coisometry belonging to $Hom(_M(H_\mu)_M, _M(H_\phi)_M)$.
Thus $[\rho_\mu]$ contains $[\rho_\phi]$. Q.E.D.

In [L3], the second author proved that for an arbitrary infinite factor $M$
(with separable pre-dual), there exists an injective subfactor $R \subset M$
satisfying $R' \cap J_M R' J_M={\bf C}$.
A subfactor $R$ of $M$ is called simple if
$R' \cap J_M R' J_M={\bf C}$.
A simple subfactor
$R$ determines the automorphisms of $M$ in the following sense;
if $\alpha, \beta \in Aut(M)$ satisfying $\alpha|_R=\beta|_R$,
then $\alpha=\beta$.
Indeed, let $u$ be the canonical implementation of $\alpha^{-1} \cdot \beta$.
Then $u\in R'$, and $u$ commutes with $J_M$.
So $u$ is a scalar, that means $\alpha=\beta$.
We can generalize this to some class of endomorphisms as follows:

\proclaim Proposition 2.10.
Let $M$ be an infinite factor and $R$ a simple  subfactor.
For every $\rho \in \En M$ with $E\in \ce M{\rho(M)}$, the following holds:
$$\{T\in M;Tx=\rho(x)T, \enskip x\in R\}=(id, \rho).\eqno (2.2)$$\par

\noindent{\it Proof.} First, we show that the general case can be reduced
to the case where $(id, \rho)=\{0\}$.
Indeed, let $\{V_i\}_i$ be an orthonormal basis of $(id, \rho)$, and $W$
an isometry in $M$ satisfying $WW^*=1-\sum V_iV_i^*$.
Then $\rho(x)=\sum V_ixV_i^*+W\sigma(x)W^*$, where $\sigma \in \En M$ is
defined by $\sigma(x)=W^*\rho(x)W$.
Note that $(id, \sigma)=\{0\}$ by construction.
If $T$ is in the left-hand side of (2.2), then $c_i:=V_i^*T\in R'\cap M=
{\bf C}$, and $W^*T$ satisfies $W^*Tx=\sigma(x)W^*T$, $x\in R$.
Since $T=\sum V_iV_i^*T+WW^*T$, if the statement is true for $\sigma$, i.e.
$W^*T=0$, we get $T=\sum c_iV_i \in (id, \rho)$.

Secondly, we construct the ``canonical implementation'' of $\rho$ as follows.
Let $\Omega$ be a separating and cyclic vector for $M$, and $L^2(M, \Omega)_+$
the natural cone with respect to $\Omega$.
Then there are unique vectors $\xi_0, \xi_1 \in L^2(M, \Omega)_+$ satisfying
$$<E(x)\Omega|\Omega>=<x\xi_0|\xi_0>.$$
$$<\rho(x)\Omega|\Omega>=<x\xi_1|\xi_1>.$$
Note that $\xi_0, \xi_1$ are cyclic because they belong to the natural cone
and implement faithful states.
So we can define an isometry $V_\rho$ by $V_\rho x \xi_1=\rho(x)\xi_0$.
We set $e_\rho=V_\rho V_\rho^*$, which is the Jones projection of $E$.
$V_\rho$ satisfies $V_\rho x=\rho(x)V_\rho$ and $J_M V_\rho =V_\rho J_M$.
Indeed, the first equality is obvious.
By identifying $e_\rho L^2(M)$ with $L^2(\rho(M), \xi_0)$, we get
$J_{\rho(M)}V_\rho=V_\rho J_M$.
On the other hand, since $e_\rho$ is the Jones projection, we have
$e_\rho J_M=J_Me_\rho =J_{\rho(M)}$.
So $V_\rho$ commutes with $J_M$.

Now suppose that $(id, \rho)=\{0\}$ and there exists a non-zero element $T$
in the left-hand side of (2.2).
Since $T^*T\in M\cap R'={\bf C}$, we may assume that $T$ is an isometry.
We set $\tilde{T}=TJ_M TJ_M$, which commutes with $J_M$ and satisfies
$\tilde{T}x=\rho(x)\tilde{T}$, $x\in R$.
Then $V_\rho^* \tilde{T} \in R'\cap J_M R'J_M={\bf C}$.
Let $\lambda=V_\rho^*\tilde{T}$, which is not zero because
$$\eqalign{<V_\rho^*\tilde{T}\xi_0|\xi_1>
&=<\tilde{T}\xi_0|\xi_0>=<T\xi_0|J_MT^*\xi_0>\cr
&=<T\xi_0|\De_\vp^{1/2}T\xi_0>=||\De_\vp^{1/4}T\xi_0||^2,}$$
where $\vp(x)=<E(x)\Omega|\Omega>$, $x\in M$.
We define a unital completely positive map $\phi :M\longrightarrow M$
by $\phi(x)=T^*\rho(x)T$, $x\in M$, which equals to
$\tilde{T}^*\rho(x)\tilde{T}$.
By construction, $[\rho]$ contains $[\rho_\phi]$.
So we show that $[\rho_\phi]$ contains $[id]$ and get contradiction.
Thanks to Proposition 2.9, it suffices to show that $\phi -|\lambda|^2id$
is completely positive.
In fact,
$$\eqalign{\phi(x)&=\tilde{T}^*\rho(x)\tilde{T}
=\tilde{T}^*e_\rho\rho(x)\tilde{T} +\tilde{T}^*(1-e_\rho)\rho(x)\tilde{T}\cr
&=\tilde{T}^*V_\rho V_\rho^*\rho(x)\tilde{T}
+\tilde{T}^*(1-e_\rho)\rho(x)\tilde{T}
=\tilde{T}^*V_\rho x V_\rho^*\tilde{T}
+\tilde{T}^*(1-e_\rho)\rho(x)\tilde{T} \cr
&=|\lambda|^2  x +\tilde{T}^*(1-e_\rho)\rho(x)\tilde{T}.} $$
Since $e_\rho$ commutes with $\rho(M)$,
$x\mapsto \tilde{T}^*(1-e_\rho)\rho(x)\tilde{T}$ is a complete positive map.
So $[\rho]$ contains $[id]$ and we get contradiction. Q.E.D.

\proclaim Corollary 2.11.
Let $M, R, \rho$ be as above, and $\sigma \in \En M$ with $d(\sigma)<\infty$.
Then the following hold:
\item {(i)} $\{T\in M;T\sigma(x)=\rho(x)T, \enskip x\in R\}=(\sigma,\rho)$.
\item {(ii)} If $\sigma|_R=\rho|_R$, then $\sigma =\rho$.

\noindent{\it Proof.} (i): Let $T$ be in the left-hand side of (i), and set
$X=\overline{\sigma}(V)R_\sigma$, where $R_\sigma$ is the isometry in (2.1).
Then $X$ satisfies $Xx=\overline{\sigma}\cdot \rho(x)X$, $x\in R$.
So thanks to Proposition 2.10, we get
$X\in (id, \overline{\sigma}\cdot \rho)$.
By simple computation using (2.1), we obtain
$V=d(\sigma)\bR \sigma^*\sigma(X)$, and $V\in (\sigma, \rho)$.
(ii): Thanks to (i), $1\in (\sigma, \rho)$, that means
$\sigma=\rho$. Q.E.D.

Let $\psi$ be a dominant weight on $M$ [CT].
Since every dominant weight is unitary equivalent,
for every $\alpha \in Aut(M)$ there is a unitary $u\in M$ satisfying
$\psi \cdot \alpha \cdot Ad(u)=\psi$.
This fact is used to define the Connes-Takesaki module of $\alpha$.
The endomorphism version is given as follows, which will be used in
the next section.

\proclaim Lemma 2.12.
Let $M$ be an infinite factor.
Then the following hold.
\item {(i)} For every $\rho \in \En M$ with $d(\rho)<\infty$, there exist a
dominant weight $\psi_\rho$ and a unitary $u\in M$ such that
$$\psi_\rho \cdot \rho\cdot Ad(u)=d(\rho)\psi_\rho, \quad
\psi_\rho \cdot E_\rho=\psi_\rho.$$
\item {(ii)} Let $\psi$ be a dominant weight.
Then for every $[\rho]\in \Se M$ with $d(\rho)<\infty$, there exists
a representative $\rho$ satisfying
$$\psi \cdot \rho=d(\rho)\psi, \quad \psi \cdot E_\rho=\psi.$$ \par

\noindent{\it Proof.} (i): Let $\psi_0$ be a dominant weight on $\rho(M)$.
Since both $d(\rho)\psi_0\cdot E_\rho$ and $\psi_0\cdot \rho$ are
dominant weights on $M$, there exists a unitary $u\in M$ satisfying
$d(\rho)\psi_0\cdot E_\rho=\psi_0\cdot \rho\cdot Ad(u)$.
So $\psi_\rho:=\psi_0\cdot E_\rho$ is the desired weight.
(ii) follows from (i) and the fact that every dominant weight is unitary
equivalent. Q.E.D.

\vskip15pt
\noindent{\bf 3. Galois Correspondence.}

In this section, we investigate the structure of irreducible
inclusions of factors with normal conditional expectations.
We present the ultimate form of the Galois correspondence of outer actions
of discrete groups and minimal actions of compact groups on factors, which
has been studied by several authors [AHKT][Ch][K][N][NT].
The key argument is how to show the existence of a conditional expectation
for every intermediate subfactor.

Let $M\supset N$ be an irreducible inclusion, i.e. $M\cap N'={\bf C}$, of
infinite factors with a conditional expectation $E\in \ce MN$.
For $\rho \in \En N$, we set
$$\cHr=\{V\in M;Vx=\rho(x)V, \enskip x\in N\}.$$
Then thanks to the irreducibility of $M\supset N$, $\cHr$ is a Hilbert space
with inner product $<V|W>1=W^*V$ as usual.
We denote by $s(\cHr)$ the support of $\cHr$, that is $\sum_i V_iV_i^*$
where $\{V_i\}_i$ is an orthonormal basis of $\cHr$.
Let $M_1$ be the basic extension of $M$ by $N$, and $e_N$ the Jones
projection of $E$.
Then $\cHr^*e_N\cHr \subset M_1\cap N'$.

Let $\gamma :M \longrightarrow N$ be the canonical endomorphism [L1][L2][L3].
Then it is known that $_NL^2(M)_N\simeq _N(_{\gamma|_N}L^2(N))_N$.
When $Ind \;E<\infty$, it is easy to show that an irreducible sector
$[\rho]\in \Se N$ is contained in $[\gamma|_N]$ if and only if $\cHr\neq 0$
(Frobenius reciprocity) [I2].
First, we establish the infinite index version of this statement.
For this purpose, it is convenient to give explicit correspondence between
submodules of $_NL^2(M)_N$ and sub-sectors of $\gamma|_N$.
Let $p\in M_1\cap N'$ be a non-zero projection.
Since both $e_N$ and $p$ are infinite projection in $M_1$, there is a partial
isometry $W\in M_1$ satisfying $WW^*=e_N$, $W^*W=p$.
Due to $e_NM_1e_N=e_NN$, we can define $\rho \in \En N$ by $WxW^*=e_N\rho(x)$,
$x\in N$.

\proclaim Lemma 3.1. Under the above assumption and notation, the following
holds:
$$_N(pL^2(M))_N\simeq _N(_\rho L^2(N))_N.$$ \par

\noindent{\it Proof.} We regard $W$ as a surjective isometry from $pL^2(M)$
to $e_NL^2(M)=L^2(N)$.
Since $M_1=J_MN'J_M$, $W$ commutes with $J_M N J_M$.
So for $\xi \in pL^2(M)$, $x,y\in N$, we obtain
$$W(x\cdot \xi \cdot y)=WxJ_My^*J_M\xi=\rho(x)WJ_My^*J_M\xi
=\rho(x)J_My^*J_MW\xi. $$
By using $e_NJ_M=J_Me_N=J_N$,
we get $W(x\cdot \xi \cdot y)=\rho(x)J_Ny^*J_NW\xi$. Q.E.D.

\proclaim Proposition 3.2. Let $M\supset N$ be an irreducible inclusion of
infinite factors with $E\in \ce MN$, and $\gamma:M \longrightarrow N$ the
canonical endomorphism.
Then for $\rho \in \En M$, the following two statements are equivalent:
\item {(i)} $\cHr\neq 0$ and the support of $E(s(\cHr))$ is 1.
\item {(ii)} $\ce N{\rho(N)}$ is non-empty and $[\rho]$ is contained in
$[\gamma|_N]$ up to multiplicity, i.e. there is decomposition
$[\rho]=\oplus [\rho_a]$ such that each $[\rho_a]$ is contained in
$[\gamma|_N]$. \par

\noindent{\it Proof.} (i) $\Rightarrow$ (ii):
Assume that $\rho$ satisfies (i).
By a simple argument, one can show that there is decomposition
$[\rho]=\oplus [\rho_a]$ such that for every $a$ there exists
$V_a \in \cH_{\rho_a}$ satisfying $E(V_aV_a^*)\geq 1$.
We set $W_a=e_N E(V_aV_a^*)^{-1/2}V_a$.
Then $W_a$ satisfies $W_aW_a^*=e_N$, $p_a:=W_a^*W_a\in M_1\cap N'$.
Since $W_axW_a^*=e_N\rho_a(x)$, $x\in N$, $[\rho_a]$ is contained in
$[\gamma|_N]$.
$\dow(p_a)=V_a^*E(V_aV_a^*)^{-1}V_a<\infty$ implies $\ce {p_aM_1p_a}{p_aN}\neq
\emptyset$, and consequently $\ce N{\rho_a(N)} \neq \emptyset$.
(ii)$\Rightarrow$ (i): It is easy to show that if $[\rho]=\oplus [\rho_a]$ and
each $\rho_a$ satisfies (i), then so does $\rho$.
Assume that $[\rho]$ is contained in $[\gamma|_N]$ and $\ce N{\rho(N)} \neq
\emptyset$.
Let $p\in M_1 \cap N'$ be the projection corresponding to $[\rho]$.
Then $\ce {pM_1p}{pN} \neq \emptyset$, which implies $p\in A\oplus B$ where
$A$ and $B$ are as in Lemma 2.7.
Let $z$ be the central support of $p$ in $M_1\cap N'$.
Since $\ma t{E\cdot \dow}$ is trivial on the center of $A\oplus B$,
$\dow|_{z(M_1\cap N')}$ is semifinite.
So there are two families of projections $\{p_a\}$, $\{q_a\}$ in
$z(M_1\cap N')$ such that $p=\sum_a p_a$, $p_a \sim q_a$ in $z(M_1\cap N')$
and $q_a \in m_{\dow}$.
Let $W_a$ be a partial isometry satisfying $W_aW_a^*=e_N$, $W_a^*W_a=q_a$,
and $\rho_a\in \En N$ defined by $W_axW_a^*=e_N\rho(x)$, $x\in N$.
Then $[\rho]=\oplus[\rho_a]$.
Since $W_a=e_NW_aq_a \in m_{\dow}$, due to Lemma 2.2, there exists $V_a\in M$
satisfying $W_a=e_NV_a$.
It is easy to check $V_a\in \cH_{\rho_a}$ and $E(V_aV_a^*)=1$.
So $\rho_a$ satisfies (i). Q.E.D.

Let $\{[\rx]\}_{\xi \in \Xi}$ be the set of irreducible sectors with finite
dimension contained in $[\gamma|_N]$.
We arrange the index set $\Xi$ such that $\overline{[\rx]}=[\rbx]$ holds for
every $\xi \in \Xi$.
For simplicity, we use notations $\Rx, \bRx, \Hx, \dx, \Ex$ instead of
$R_{\rx}, \bR {\rx},$ etc.
We define the Frobenius maps $\cx:\Hx \longrightarrow \Hbx $,
$\bcx:\Hbx \longrightarrow \Hx$ by
$$\cx(V)=\sqrt{\dx}V^*\bRx, \quad V\in \Hx,$$
$$\bcx(\bV)=\sqrt{\dx}\bV^* \Rx, \quad \bV \in \Hbx.$$
Then thanks to (2.1), $\bcx\cx=1_{\Hx}$, $\cx\bcx=1_{\Hbx}$.
So in particular, both $\cx$ and $\bcx$ are invertible.
We introduce a new inner product to $\Hx$ by
$$(V_1,V_2)1=\dx E(V_1V_2^*)\in (\rx,\rx)={\bf C}, \quad V_1, V_2 \in \Hx.$$
Due to the estimate $|(V_1,V_2)|\leq \dx ||V_1||||V_2||$, there is a
non-singular positive operator $\ax\in B(\Hx)$ satisfying
$$(V_1,V_2)=<\ax V_1|V_2>.$$
Let $\{V_i\}_i \subset \Hx$ be an orthonormal basis of $\Hx$.
Since $\sum V_iV_i^*=s(\Hx)\leq 1$, we get
$$Tr(\ax)=\sum (V_i,V_i)=\dx E(s(\Hx))\leq \dx.$$
So $\ax$ is a trace class operator.
By simple computation one can show the following:
$$<\cx(V_1)|\cx(V_2)>=(V_1,V_2)=<\ax V_1|V_2>, $$
$$<\bcx(\bV_1)|\bcx(\bV_2)>=(\bV_1,\bV_2)=<\bax \bV_1|\bV_2>. $$
Thus we get $\cx^*\cx=\ax$, $\bcx^*\bcx=\bax$.
This shows that $\ax$ is an invertible trace class operator, that implies
$\nx:=dim \Hx<\infty$.
Thanks to $\bcx=\cx^{-1}$, we obtain
$$Tr(\bax)=Tr(\bcx^*\bcx)=Tr(\bcx\bcx^*)=Tr(\ax^{-1}).$$
This implies
$${1\over \dx}\leq \ax \leq \dx, \quad \nx \leq \dx^2.$$
If $\ax=1$, (this is the case if for instance $Ind\; E < \infty$),
then $\nx\leq \dx$.
On the other hand if $\nx=\dx$, then it is easy to show $\ax=1$ and
$E(s(\Hx))=1$, i.e. $s(\Hx)=1$.

\proclaim Theorem 3.3.
Let $M\supset N$ be an irreducible inclusion of infinite factors with
$E\in \ce MN$, and $\hrc =A\oplus B_1\oplus B_2 \oplus C$ the decomposition
described in Proposition 2.8.
Then with the same notation as above, the following hold:
\item {(i)} $A=\oplus_{\xi \in \Xi} \Ax$, where $\Ax=\Hx^* e_N \Hx \simeq
M(n_\xi, {\bf C})$.
\item {(ii)} $B_1$ and $B_2$ are of type I.
\item {(iii)} For $V_1, V_2\in \Hx$,
$\ma t{E\circ \dow}(V_1^*e_NV_2)=V_1^*\ax^{-it}e_N \ax^{it}V_2.$
\item {(iv)} For $V_1, V_2\in \Hx$,
$j(V_1^*e_NV_2)=\cx(\ax^{1/2}V_2)^*e_N \cx(\ax^{1/2}V_1).$
\par

\noindent {\it Proof.}
(i): Thanks to Proposition 2.8, $A$ is direct sum of type I factors.
By using the one-to-one correspondence as described just before Lemma 3.1,
we can parametrize the direct summands of $A$ by $\Xi$ such that
$A=\oplus A_\xi$ and $\Ax\supset \Hx^* e_N \Hx$ hold.
So it suffices to show that $\Ax$ is of type I$_{n_\xi}$.
If $\Ax$ is finite, then $\Ax \subset m_{\dow}$ because $\dow|_{\Ax}$ is
semifinite.
So we can take matrix units $\{e_{i,j}\}_{1\leq i,j}$ of $\Ax$
(with $\sum e_{i,i}=1_{\Ax}$) such that $\dow(e_{i,j})=b_i\delta_{i,j}$.
We may assume that there is a partial isometry $W_1\in M_1$ satisfying
$W_1W_1^*=e_N$, $W_1^*W_1=e_{1,1}$ and $W_1xW_1^*=e_N\rx(x)$ for $x\in N$.
We set $W_i=W_1e_{1,i}$.
Then thanks to Lemma 2.2, there exists $V_i\in M$ such that
$W_i=\sqrt{b_i}e_N V_i$.
$\{V_i\}$ is an orthonormal basis of $\Hx$.
Indeed, it is easy to show that it is an orthonormal system.
Suppose $V\in \Hx$ is perpendicular to $\{V_i\}$.
Since $e_{i,j}=W_i^*e_NW_j=\sqrt{b_ib_j}V_i^*e_NV_j$, $V^*e_NV$ is an
element in $\Ax$ satisfying $e_{i,j}V^*e_NV=0$.
This means $V^*e_NV=0$ and $0=\dow(V^*e_NV)=V^*V$, i.e. $V=0$.
So $\{V_i\}$ is an orthonormal basis of $\Hx$ and the rank of $\Ax$ coincides
with $n_\xi$.
Now suppose $\Ax$ is of type I$_\infty$.
Since $\dow|_{\Ax}$ is semifinite, there is a matrix unit
$\{e_{i,j}\}_{1\leq i,j <\infty}$ (not necessarily $\sum e_{i,i}=1$),
such that $\dow(e_{i,i})<\infty$, $\dow(e_{i,j})=0$ for $i\neq j$.
Then we can define $W_i$ and $V_i$ as before.
However, $\{V_i\}_{1\leq i<\infty}$ is an orthonormal system of $\Hx$,
that contradicts the fact $dim \Hx=n_\xi<\infty$.

\noindent (ii):
Since $\dow|_{B_1}$ is semifinite and $j(B_1)=B_2$, it suffices to show
that $pB_1p$ is of type I for every $p\in B_1$ with $\dow(p)<\infty$.
Let $W \in M_1$ be a partial isometry with $WW^*=e_N$, $W^*W=p$, and
define $\rho\in \En M$ by $WxW^*=e_N\rho(x)$, $x\in N$ as before.
Thanks to Lemma 2.2, there exists an isometry $V\in \Hr$ satisfying
$W=\sqrt{c}e_NV$, $c=\dow(p)$.
So $E(VV^*)={1\over c}$ and we get ${1\over c}\leq E(s(\Hr))\leq 1$.
Let $P=N\cap \rho(N)'$.
Then in the same way as in the proof of Lemma 3.1, we can show that
$pB_1p$ is isomorphic to $P$.
So we show that $P$ is of type I.
Thanks to $P\Hr=\Hr$, we can define a normal representation of $P$ on $\Hr$ by
$\pi(x)V=xV$, $x\in P$, $V\in \Hr$.
Note that ${1\over c}\leq E(s(\Hr))$ implies that $\pi$ is faithful.
Thus to prove that $P$ is of type I, we show that there exists a normal
conditional expectation from $B(\Hr)$ to $\pi(P)$ [S, Proposition 10.21].
For $\omega \in P_*$ we can define a bilinear form on $\Hr$ by
$\omega(E(V_1V_2^*))$, $V_1, V_2\in \Hr$ with an estimate
$|\omega(E(V_1V_2^*))|\leq ||\omega||||V_1||||V_2||$.
So there exists a unique bounded operator $h_\omega$ satisfying
$$\omega(E(V_1V_2^*))=<h_\omega V_1|V_2>.$$
For $x,y\in P$, $\omega\in P_*$, $h_\omega$ satisfies
$h_{x\cdot \omega\cdot y}=\pi(x)h_\omega \pi(y)$.
Indeed, by definition we get
$$\eqalign{x\cdot \omega \cdot y(E(V_1V_2^*))
&=\omega (yE(V_1V_2^*)x)=\omega(E(\pi(y)V_1(\pi(x)^*V_2)^*))\cr
&=<h_\omega\pi(y)V_1|\pi(x)^*V_2>=<\pi(x)h_\omega\pi(y)V_1|V_2>.}$$
If $\omega \in P_*$ is positive, we have
$$Tr(h_\omega)=\omega(E(s(\Hr))) \leq\omega(1)=||\omega||,$$
so by using polar decomposition of linear functionals and the fact just proved
above, we get
$$||h_\omega||_1:=Tr(|h_\omega|)=Tr(h_{|\omega|})
\leq ||(|\omega|)||=||\omega||, \quad \omega\in P_*.$$
Hence we can define a bounded order preserving linear map
$\theta :P_*\longrightarrow B(\Hr)_*$ by
$\theta(\omega)(a)=Tr(h_\omega a)$, $a\in B(\Hr)$.
Note that $\theta$ satisfies $\theta(x\cdot \omega\cdot y)
=\pi(x)\cdot \theta(\omega)\cdot \pi(y)$, $x,y\in P$.
Let $F_0$ be the transposition of $\theta$.
Then $F_0$ is a positive normal map $F_0 :B(\Hr) \longrightarrow P$
satisfying $F_0(\pi(x)a\pi(y))=xF_0(a)y$, $x,y\in P$, $a\in B(\Hr)$.
Note that $F_0(1)=E(s(\Hr))$ is a central element of $P$ because
$us(\Hr)u^*=s(\Hr)$ holds for every unitary $u\in P$.
Since $E(s(\Hr))$ is invertible, we can define a normal conditional
expectation $F:B(\Hr)\longrightarrow \pi(P)$ by
$$F(a)=\pi(E(s(\Hr))^{-1/2}F_0(a)E(s(\Hr))^{-1/2}),\quad a\in B(\Hr).$$
Therefore,  $P$ is of type I.

\noindent (iii):
By a simple argument one can show that unitary perturbation of $\rx$
does not have any effect on the formulae in (iii)(iv).
So thanks to Lemma 2.12, we assume that there is a dominant weight $\psi$ on
$N$ satisfying $\psi \cdot \rx =d(\xi)\psi$, $\psi \cdot E_\xi=\psi$
for every $\xi \in \Xi$.
Then $\ma t\psi$ commute with $\rx$ and we get $\ma t{\psi\cdot E}(\Hx)=\Hx$.
So we show $\ma t{\psi\cdot E}(V)=\ax^{it}V$ for $V\in \Hx$, that implies
the statement.
Indeed, since $dim \Hx<\infty$, every element in $\Hx$ is analytic for
$\{\ma t{\psi\cdot E}\}$.
Let $V\in \Hx$ and $x\in m_\psi$.
Then by using the KMS condition, we obtain
$$\psi \cdot E(VxV^*)=\psi \cdot E(xV^*\ma {-i}{\psi \cdot E}(V))
=<\ma {-i}{\psi \cdot E}(V)|V>\psi(x).$$
On the other hand, from $E(VxV^*)=E(\rx(x)VV^*)={1\over d(\xi)}(V,V)\rx(x)$
we get
$$\psi\cdot E(VxV^*)={1\over d(\xi)}(V,V)\psi\cdot \rx(x)=<\ax V|V>\psi(x).$$
So we obtain $\ma t{\psi\cdot E}(V)=\ax^{it}V$.

\noindent (iv):
Let $\zx$ be the unit of $\Ax$.
Then by using the correspondence between sub-bimodules of $_NL^2(M)_N$ and
sub-sectors of $\gamma_N$, we get $j(\Ax)=\Abx$ and $j(\zx)=\zbx$.
Let $\psi$ be as before.
Then due to (i), it is easy to show that
$\Hx^*\Lambda_{\psi \cdot E}(n_\psi \cap n_\psi^*)$ is dense in
$\zx H_{\psi \cdot E}$.
Since both $j(V_1^*e_NV_2)$ and $\cx(\ax^{1/2}V_2)^*e_N \cx(\ax^{1/2}V_1)$
belong to $\Abx$, it suffices to show the equality on
$\Hbx^*\Lambda_{\psi \cdot E}(n_\psi \cap n_\psi^*)$.
Let $a\in n_\psi \cap n_\psi^*$ and $X\in \Hbx$.
Since $V_1, V_2$ are analytic elements for $\{\ma t{\psi\cdot E}\}$, we get
$$\eqalign{j(V_1^*e_NV_2)\Lambda_{\psi \cdot E}(X^*a)
&=J_MV_2^*J_Me_NJ_MV_1J_M \Lambda_{\psi \cdot E}(X^*a)\cr
&=J_MV_2^*J_Me_N\Lambda_{\psi \cdot E}
(X^*a\ma {i\over 2}{\psi \cdot E}(V_1)^*)\cr
&=J_MV_2^*J_Me_N\Lambda_{\psi \cdot E}
(X^*\ma {i\over 2}{\psi \cdot E}(V_1)^*\rx(a))\cr
&=J_MV_2^*J_M \Lambda_{\psi \cdot E}
(E(X^*\ma {i\over 2}{\psi \cdot E}(V_1)^*)\rx(a))\cr
&=\Lambda_{\psi \cdot E}(E(X^*\ma {i\over 2}{\psi \cdot E}(V_1)^*)\rx(a)
\ma {-{i\over 2}}{\psi \cdot E}(V_2))\cr
&=\Lambda_{\psi \cdot E}(E(X^*\ma {i\over 2}{\psi \cdot E}(V_1)^*)
\ma {-{i\over 2}}{\psi \cdot E}(V_2)a).}$$
By using $X=\cx(\bcx(X))=\sqrt{d(\xi)}\bcx(X)^*\bRx$, we get
$$\eqalign{j(V_1^*e_NV_2)\Lambda_{\psi \cdot E}(X^*a)
&=\sqrt{d(\xi)}\Lambda_{\psi \cdot E}
(\bRx^*E(\bcx(X)\ma {i\over 2}{\psi \cdot E}(V_1)^*)
\ma {-{i\over 2}}{\psi \cdot E}(V_2)a)\cr
&={1\over \sqrt{d(\xi)}}(\bcx(X), \ma {i\over 2}{\psi \cdot E}(V_1))
\Lambda_{\psi \cdot E}(\bRx^* \ma {-{i\over 2}}{\psi \cdot E}(V_2)a)\cr
&={1\over d(\xi)}(\bcx(X), \ax^{-1/2}V_1)
\Lambda_{\psi \cdot E}(\cx(\ax^{1/2}V_2)^*a).}$$
On the other hand, we have
$$\cx(\ax^{1/2}V_2)^*e_N\cx(\ax^{1/2}V_1)\Lambda_{\psi \cdot
E}(X^*a)={1\over d(\xi)}(\cx(\ax^{1/2}V_1), X)
\Lambda_{\psi \cdot E}(\cx(\ax^{1/2}V_2)^*a),$$
so it suffices to show $(\bcx(X), \ax^{-1/2}V_1)=(\cx(\ax^{1/2}V_1), X)$.
Actually,
$$(\cx(\ax^{1/2}V_1), X)=<\bcx(X)|\ax^{1/2}V_1>=(\bcx(X), \ax^{-1/2}V_1).$$
Q.E.D.

\noindent {\it Remark 3.4.} Let $V_1, V_2 \in \Hx$.
Then we get
$$\dow(V_1^*e_NV_2)=<V_2|V_1>.$$
$$\dow(j(V_1^*e_NV_2))=<\cx(\ax^{1/2}V_1)|\cx(\ax^{1/2}V_2)>=
(\ax^{1/2}V_2,\ax^{1/2}V_1)=<\ax^2V_2|V_1>.$$
So $\dow\cdot j|_A=\dow|_A$ if and only if $\ax=1$ for all $\xi\in \Xi$.
It is also easy to show that $\dow|_A$ is a trace if and only if $\ax$ is a
scalar for all $\xi\in \Xi$.

To the best knowledge of the authors there is no known example which violates
$\ax=1$.
However, the following example shows that $\dow \cdot j|_{\hrc}=\dow|_{\hrc}$
does not hold in general; $B_i$, $i=1,2$ may not vanish.

\proclaim Example 3.5. \par
\noindent (i) Let $G$ be a discrete group and $H$ a subgroup, and let
$\alpha$ be an outer action of $G$ on a factor $L$.
We set $N=L\times_\alpha H$, $M=L\times_\alpha H$.
Then $M\supset N$ is an irreducible inclusion of factors with a unique
conditional expectation $E$.
We identify $M$ and $N$ with $(L\otimes{\bf C})\times G$ and
$(L\otimes {\bf C})\times H$ acting on
$L^2(L)\otimes \ell^2(G/H)\otimes \ell^2 (G)$ in an obvious sense.
Let $f$ be the orthogonal projection onto ${\bf C}\delta_{\dot e}\subset
\ell^2(G/H)$, where $\delta$ stands for the $\delta$-function
and $\dot{e}$ the class of the neutral element $e$,
and $F_0=id \otimes Tr\in \ow {L\otimes \ell^2(G/H), L\otimes {\bf C}}$.
Then thanks to Lemma 2.3 we can identify $M_1$ with
$(L\otimes \ell^\infty (G/H))\times G$ where the action of $G$ on
$\ell^\infty(G/H)$ is the translation,
$e_N$ with $1\otimes f \otimes 1$ and $\dow$ with the natural
extension of $F_0$ to $(L\otimes\ell^\infty(G/H)) \times G$.
So under this identification we get $\hrc=\ell^\infty(H\backslash G/H)$.
For $\dot{g}\in G/H$, we denote by
$p_{\dot{g}}\in \ell^\infty(H\backslash G/H)$ the
projection corresponding to the $H$-orbit of $\dot{g}$.
Then $\dow(p_{\dot{g}})$ is exactly the length of the orbit, i.e.
$\dow(p_{\dot{g}})=[H:H_g]$ where $H_g:=gHg^{-1}\cap H$.
$j(p_{\dot{g}})$ can be computed by using bimodules as in [KY], and we have
$j(p_{\dot{g}})=p_{\dot{g^{-1}}}$.
So for example if $g^{-1}Hg\subset H$ and $g^{-1}Hg\neq H$, then
$\dow(p_{\dot{g}})=1$ although $\dow(j(p_{\dot{g}}))\neq 1$.
Let $G$ be the group generated by the finite permutations of
${\bf Z}$ and $g$ where $g$ is the translation of ${\bf Z}$, and $H$ the
finite permutations of ${\bf N}\cup \{0\}$.
Then $gHg^{-1}$ is the finite permutation of {\bf N} and we get
$gHg^{-1}\subset H$, $[H:H_g]=\infty$.
So we obtain $\dow(p_{\dot{g}})=\infty$, $\dow(p_{\dot{g^{-1}}})=1$.
This means $B_i\neq \{0\}$, $i=1,2$ in this example.

\noindent (ii) Let $G\supset H$ be a pair of discrete groups with the
following property: for evry $g\neq e\in G$ $\{hgh^{-1};h\in H\}$ is an
infinite set.
Let $M:=L(G)$ be the group von Neumann algebra of $G$ and $N:=L(H)$ the
subfactor of $M$ generated by $H$.
Then in exactly the same way as one proves that $M$ is a factor, one can show
$M\cap N'={\bf C}$.
Although this example looks similar to the previous one, these two have
essencially different natures.
As before we can identify $N$, $M$ and $M_1$ with ${\bf C}\times H$,
${\bf C}\times G$ and $\ell^\infty(G/H)\times G$ acting on
$\ell^2(G/H)\otimes \ell^2(G)$.
However, we can conclude only $\ell^\infty(H\backslash G/H)\subset \hrc$
because the action of $G$ on $G/H$ is not necessarily free.
In fact the equlity does not hold in general.
For example, let $G=F_3$ be the free group generated by $g_1,g_2,g_3$ and
$H=F_2=<g_1, g_2>$.
Then the $N-N$ bimodule $_NX_N$ generated by $\delta_{g_3} \in \ell^2(F_3)$
is equivalent to $_N \ell^2(F_2)\otimes \ell^2(F_2)_N$ where $\otimes$ is
the usual tensor product and the left and the right actions act on each tensor
component respectively.
So $End(_NX_N)\simeq N^{op}\otimes N$.
This means that $\hrc$ has a type II summand.
Actually a little more effort shows that
$A={\bf C}e_N$, $B_1=B_2={0}$ and $C$ is of type II where $A$, $B_1$, $B_2$,
and $C$ are as in Proposition 2.8.

\noindent {\it Remark 3.6.}
Let $M_2$ be the basic extension of $M_1$ by $M$ in the first exmple,
i.e. $M_2:=J_{M_1}M'J_{M_1}$.
Then it is easy to show
$M_2=L\otimes B(\ell^2(G/H))\times_{\alpha\otimes Ad(\pi)} G$
where $\pi$ is the translation.
So $M'\cap M_2=\pi(G)'$.
In [B], W. Binder constructs an example of a pair of discrete groups
$G \supset H$ such that $\pi(G)'$ is a type III factor.
This means that the restriction of the unique expectation in $\ce {M_2}{M_1}$
to $M'\cap M_2$ may fail to be a trace in general.

In [HO] R. Herman and A. Ocneanu called an inclusion of factors $M\supset N$
discrete if $\ce MN$ is not empty.
However, the above examples show that existence of a normal conditional
expectation is not strong enough to assure properties resembling those of
crossed products by discrete group actions.
Therefore, in this paper we use the terminology in the following sense.

\proclaim Definition 3.7. An inclusion of factors is called
discrete if $\ce MN$ is non-empty and $\dow|_{\hrc}$ is semifinite
for some  (and equivalently all) $E\in  \ce MN$. \par

In what follows we assume that $M\supset N$ is an irreducible discrete
inclusion of infinite factors.
Note that discreteness is equivalent to $\hrc =A$ in the decomposition given
in Proposition 2.8, and to $[\gamma|_N]=\oplus \nx [\rx]$, $d(\xi)<\infty$.

For each $\xi \in \Xi$ choose an orthogonal basis $\{\Vxi\}_{i=1}^{\nx}$
consisting of eigenvectors of $\ax$ belonging to $\axi$.
For $x\in M$ we define the ``Fourier coefficient'' $\xxi$ by
$$\xxi={d(\xi)\over \axi} E(\Vxi x).$$
Then $x$ has the following formal expansion:
$$x =\sum_{\xi\in \Xi}\sum_{i=1}^{\nx}\Vxi^*\xxi.$$
Although the above sum does not converge even in weak topology in general,
we can give justification of the expansion as follows.
We define $\pxi \in \hrc$ by
$$\pxi={d(\xi)\over \axi}\Vxi^* e_N \Vxi.$$
Then $\pxi$ is a projection with $z_\xi=\sum_{i=1}^{\nx}\pxi$,
where $\zx$ is the unit of $\Ax$.
By discreteness assumption we have $\sum_{\xi \in \Xi}\zx =1$.
Let $\omega$ be a faithful normal state on $N$ and set $\vp=\omega \cdot E$.
Since $\pxi\La_\vp(x)=\La_\vp(\Vxi^* \xxi)$ and
$\La_\vp(x)=\sum_{\xi, i} \pxi \La_\vp(x)$, the sum converges in Hilbert
space topology.
Note that $\{x(\xi)_i\}$ uniquely determines $x$ while it is difficult to tell
when a series $\{x(\xi)_i\}$ is actually the Fourier coefficient of some
element $x\in M$.

Although the following lemma might sound trivial, we need to prove it
because the expansion does not make sense in any decent operator algebra
topology.

\proclaim Lemma 3.8. Under the above assumption, assume that there is an
assignment of subspaces $\Kx \subset \Hx$ satisfying the following conditions.
\item {(i)} $\ax\Kx\subset \Kx.$
\item {(ii)} $\Kx^*\subset N{\cal K}_{\bar \xi}.$
\item {(iii)} Let $\eta,\zeta \in \Xi$ and set
$\Xi_{\eta,\zeta}=\{\xi \in \Xi;\rho_\eta \rho_\zeta \succ \rx\}.$
Then, $\cK_\eta\cK_\zeta \subset \sum_{\xi \in \Xi_{\eta,\zeta}}N\Kx.$

\noindent{\it
Let $L$ be the von Neumann algebra generated by $N$ and
$\{\Kx\}_{\xi \in \Xi}$.
Then there exists $E_L\in \ce ML$, and $L$ is characterized by
$$L=\{x\in M; E(\Kxo x)=0, \quad \xi \in \Xi\},$$
where $\Kxo$ is the orthogonal complement of $\Kx$ with respect to $<|>$. }

\noindent {\it Proof.}
Let $L_0$ be the direct sum of $\Kx^*N$.
Thanks to (ii) and (iii), $L_0$ is the *-algebra generated by $N$ and
$\{\Kx\}$, which is dense in $L$.
Let $L_1=\{x\in M; E(\Kxo x)=0, \quad \xi \in \Xi\}$ and $K$ the closure of
$\La_\vp(L)$ in $H_\vp$.
First, we claims $L_1=\{x\in M;\La_\vp(x)\in K\}.$
Indeed, due to (i) we may arrange $\{\Vxi\}$ such that $\{\Vxi\}_{i=1}^{\mx}$
is an orthonormal basis of $\Kx$.
Then we get $K=\oplus_{\xi \in \Xi}\oplus_{i=1}^{\mx} H_{\xi,i}$ where
$H_{\xi,i}=\pxi H_\vp$, and so
$$L_1=\{x\in M; \pxi\La_\vp(x)=0, \quad i>\mx\}.$$
Thus we get the claim.
Secondly, we show that there exists $E_L\in \ce ML$ with $\vp \cdot E_L=\vp$.
Thanks to the Takesaki theorem on conditional expectations [S], it suffices to
prove $\ma t\vp(L)=L$, or in our case $\ma t\vp(\Kx)\subset N\Kx$.
As before we may and do assume that there is a dominant weight $\psi$ on $N$
satisfying $\psi \cdot E_\xi=\psi$, $\psi \cdot \rx=d(\xi)\psi$,
so we have $\ma t{\psi\cdot E}(V)=\ax^{it}V$ for $V\in \Hx$.
We set $u_t^\xi=[D\omega:D\psi]_t\rx([D\omega:D\psi]_t^*)\in N$, where
$[D\omega:D\psi]_t$ is the Connes cocycle derivative.
Then we get
$$\ma t\vp(V)=Ad([D\omega\cdot E:D\psi\cdot E]_t)\cdot
\ma t{\psi\cdot E}(V)
=Ad([D\omega:D\psi]_t)(\ax^{it}V)=u_t^\xi \ax^{it}V,$$
so due to (i) we get $\ma t\vp(\Kx)\subset N\Kx$.
Now let $e_L$ be the Jones projection for $E_L$, i.e.
$e_L\La_\vp(x)=\La_\vp(E_L(x))$, $x\in M$.
Then $e_L$ is the orthogonal projection onto $K$.
Since $L$ is characterized by
$L=\{x\in M; e_L\La_\vp(x)=\La_\vp(x)\},$
we get $L=L_1$. Q.E.D.

The following is the main technical result in this paper.

\proclaim Theorem 3.9. Let $M\supset N$ be an irreducible inclusion of
infinite factors with $E\in \ce MN$.
We assume that the inclusion is of discrete type and $\ma t{E\cdot \dow}$
is trivial.
Let $L$ be an intermediate subfactor and $\Kx=L\cap \Hx$.
Then $\{\Kx\}$ satisfies the assumption of Lemma 3.8 and
$L$ is generated by $N$ and $\{\Kx\}$.
Consequently, there exists $E_L\in \ce ML$.
\par

\noindent {\it Proof.} First, we show that the statement can be reduced to
the case where $N$ is of type III.
Suppose that the statement holds for type III factors.
Then we apply the statement to $\hat{M}=M\otimes P$, $\hat{N}=N\otimes P$ and
$\hat{L}=L\otimes P$ where $P$ is a type III factor,
and get that  $\hat{L}$ is generated by $\hat{N}$ and
$(\Hx\otimes {\bf C})\cap \hat{L}=\Kx\otimes {\bf C}$.
$\{\Kx\}$ satisfies the assumption of Lemma 3.8 because so does
$\{\Kx\otimes 1\}$ by assumption.
Thanks to Lemma 3.8 we get
$$\hat{L}=\{x\in M\otimes P;(E\otimes id)((\Kxo\otimes 1)x)=0,
\quad \xi \in \Xi\},$$
and so we obtain
$$L=\{x\in M; E(\Kxo x)=0, \quad \xi \in \Xi\}.$$
Therefore, the statement holds for $L$ as well.
Now, we assume that $N$ is of type III.
Let $\{\Vxi\}$ be as in the proof of Lemma 3.8.
Thanks to $\Hx^*\subset N\cH_{\bar \xi}$, $\cH_\eta \cH_\zeta\subset
\sum_{\xi\in \Xi_{\eta,\zeta}}N\Hx$ and the Fourier decomposition,
to prove that $\{\Kx\}$ satisfies the assumption of Lemma 3.8 it suffices
to show $\xxi=0$ for $x\in L$, $\xi \in \Xi$, $i>\mx$,
which is actually enough for the statement due to Lemma 3.8.
Suppose the converse; there exists $x\in L$ such that $\xxi\neq 0$
for some $\xi\in \Xi$ and some $i>\mx$.
Let $y=axb$, $a,b\in N$.
Then $E_\xi(y(\xi)_i)=\rx(a)E_\xi(\xxi b)$.
since $N$ is a type III factor, we can choose $a,b$ such that
$E_\xi(y(\xi)_i)=1$, so we assume $E_\xi(\xxi)=1$ from the beginning.
Let $R$ be a simple injective subfactor of $N$ and $U(R)$ the unitary group
of $R$.
We set $\cC=\overline{conv\{ux\rx(u^*);u\in U(R)\}}^w$ and define an action
$\theta$ of $U(R)$ on $\cC$ by $\theta_u(w)=uw\rx(u^*)$, $u\in U(R)$,
$w\in \cC$.
We claim that the set of fixed points of $\cC$ under $\theta$, which is
the same as $\{ w\in \cC;aw=w\rx(a), \quad a\in R\}$, is non-empty.
Indeed, since $R$ is AFD, there exists an increasing sequence of finite
dimensional unital von Neumann-subalgebras $\{R_n\}_{n=1}^\infty$
generating $R$.
Let $\cC_n$ be the fixed points of $\cC$ under $\theta|_{U(R_n)}$, that
is a non-empty compact set because $U(R_n)$ is a compact group.
Then $\{\cC_n\}_{n=1}^\infty$ is a decreasing sequence of non-empty compact
sets, and so $\cC_\infty:=\cap_{n=i}^\infty \cC_n$ is non-empty as well.
Let $w\in \cC_\infty$.
Then $w$ satisfies $aw=w\rx(a)$ for $a\in \cup_nR_n$, and for $a\in R$ because
$\cup_n R_n$ is dense in $R$.
Thus $\cC_\infty$ is the set of the fixed points.
From the definition of the Fourier coefficient of $w\in \cC_\infty$ we get
$\rho_\eta(a)w(\eta)_j=w(\eta)_j\rx(a)$ for $a\in R$, $\eta \in \Xi$.
Applying Corollary 2.11 we obtain $w(\eta)_j=0$ for $\eta\neq \xi$ and
$w(\xi)_j\in {\bf C}$, that means $w^*\in \Hx\cap N=\Kx$.
On the other hand,  $E_\xi(\xxi)=1$ implies $E_\xi((ux\rx(u^*))(\xi)_i)
=\rx(u)E_\xi(\xxi)\rx(u^*)=1$ for $u\in U(R)$,
and so $E_\xi(w(\xi)_i)=1$ by continuity.
Since $w(\xi)_i$ is a scalar $w(\xi)_i=1$.
Hence $w^* \notin \Kx$, that is contradiction.
Therefore we get $\xxi=0$ for $x\in L$, $\xi \in \Xi$, $i>\mx$.
Q.E.D.

\proclaim Corollary 3.10. Let $M$, $N$, $\Xi$  be as above and $\Xi_1$
a self-conjugate subset of $\Xi$ with the following properties; whenever
$\xi, \eta \in \Xi_1$, $\Xi_{\xi, \eta} \subset \Xi_1$.
Then there exists an unique intermediate subfactor $L$ such that
if we denote by $\gamma'$ the canonical endomorphism
$\gamma' : L \longrightarrow N$, then
$$[\gamma'|_N]=\oplus_{\xi \in \Xi_1}\nx [\rx].$$
\par
\noindent{\it Proof.} Set $L=N \vee \{\Hx\}_{\xi \in \Xi_1}$. Q.E.D. \par

\proclaim Corollary 3.11.
Let $M\supset N$ be an irreducible inclusion of factors ($N$ is not
necessarily infinite) with $E\in \ce MN$.
We assume that the inclusion is of discrete type and $\ma t{E\cdot \dow}$
is trivial.
Then for every intermediate subfactor $L$, $\ce ML$ is not empty. \par

\noindent{\it Proof.} It is enough to prove the statement when $N$ is finite
and $M$ is infinite.
Let $F$ be a type I$_\infty$ factor.
Then thanks to Theorem 3.9 $\ce {M\otimes F}{L\otimes F}$ is not empty.
Since we can identify $M\supset L$ with
$e(M\otimes F)e \supset e(L\otimes F)e$ where $e$ is a minimal projection of
$F$, $\ce ML$ is not empty. Q.E.D.

\noindent {\it Remark 3.12.} Using the same type of argument, we can show the
following:
for an irreducible discrete inclusion of infinite factors $M\supset N$ and
a simple injective subfactor $R$ of $N$, $M\cap R'={\bf C}$ holds.
Indeed, suppose $x\in M\cap R'$.
Then $\xxi$ satisfies $\xxi a=\rx(a)\xxi$ for $a\in R$.
So we get $\xxi=0$ unless $\rx=id$, and $x\in N\cap R'={\bf C}$.

The above theorem means that when $\ax$ is a scalar there is one-to-one
correspondence between the set of intermediate subfactors and that of
the systems of Hilbert subspaces $\{\Kx\}$ satisfying (ii) and (iii)
of Lemma 3.8.
This observation has a lot of useful applications in Galois correspondence of
operator algebras as stated below.
Although our statements can be unified as that of depth 2 irreducible
inclusions of discrete type in the language of Kac algebras
(see next section), we first state them in two classical cases:
crossed product inclusions of outer actions of discrete groups and
fixed point inclusions of minimal actions of compact groups.

\proclaim Theorem 3.13. Let $G$ be a discrete group and $\alpha$ an outer
action of $G$ on a factor $N$.
Then the map $H \mapsto N\times_\alpha H$ gives
one-to-one correspondence between the lattice of all subgroups of $G$ and
that of all intermediate subfactors of $N\subset N\times_\alpha G$. \par

\noindent{\it Proof.} Let $\{\lambda(g)\}$ denote the implementing unitaries
of $\alpha$ in $M:=N\times_\alpha G$.
Then it is easy to see $\Xi=G$ and $\cH_g= {\bf C}\lambda(g)$ where $\Xi$ and
$\cH_g$ are as in Theorem 3.9.
(Note that the argument in Theorem 3.9 makes sense even when $N$ is finite
as far as $\{\rx\}$ are automorphisms.)
Let $\{\cK_g\}_{g\in G}$ be a system of subspaces satisfying (ii) and (iii)
of Theorem 3.9.
Then there exists a subgroup $H\subset G$ such that $\cK_g={\bf C}\lambda(g)$
if $g\in H$ and $\cK_g={0}$ if $g\notin H$.
This means that for every intermediate subfactor $L$ there exits a subgroup
$H$ with $L=N\times_\alpha H$. Q.E.D.

\noindent{\it Remark 3.14.} In [Ch], H. Choda proved that there is
one-to-one  correspondence between the set of subgroups and the set
of intermediate  subfactors $L$ with $\ce {N\times_\alpha G}L$
non-empty.  The above theorem says that the existence of a normal
conditional  expectation to every intermediate subfactors
automatically follows from  Theorem 3.9.

Let $G$ be a compact group.
We call an action $\alpha$ of $G$ on a factor $M$ minimal if $\alpha$ is
faithful and $M\cap {M^G}'={\bf C}$ where $M^G$ is the fixed point algebra
under $\alpha$.
It is known that if $\alpha$ is minimal the crossed product $M\times_\alpha G$
is always a factor (See Remark 4.5).
We fix a complete system of representatives of the
equivalence classes of the irreducible representations of $G$ and denote it
by $\hG$.
If $\alpha$ is minimal and the fixed point algebra $M^G$ is infinite,
using the same type of the argument as in [AHKT, Lemma III 3.4], one can show
that for every $\pi\in \hG$ there exists a Hilbert space $\cH_\pi \in M$ with
support 1 such that $\cH_\pi$ is globally invariant under $\alpha$ and
$\alpha|_{\cH_\pi}$ is equivalent to $\pi$.
This means that $M$ is the crossed product of $M^G$ and the dual object of $G$
by the corresponding Roberts action [R1].
We fix such a $\cH_\pi$ for each $\pi \in \hG$ and choose an orthonormal
basis $\{\Vpi\}_{i=1}^{d(\pi)}$ of $\cH_\pi$ where $d(\pi)$ is the dimension
of $\cH_\pi$.
Let $N=M^G$ and $E$ the unique element in $\ce MN$ obtained by
$$E(x)=\int_G\alpha_g(x)dg, \quad x\in M.$$
We define an endomorphism $\rho_\pi\in End(N)$ by
$$\rho_\pi(x)=\sum_{i=1}^{d(\pi)}\Vpi x \Vpi^*, \quad x\in N.$$
Thanks to the minimality of $\alpha$, $\rho_\pi$ is always irreducible with
$d(\rho_\pi)=d(\pi)$.
It is routine to show that $\rho_\pi$ does not depend on the choice of the
basis and that the sector of $\rho_\pi$ does not depends on the choice of
$\cH_\pi$.
Note that $\cH_\pi$ is characterized by
$$\cH_\pi=\{V\in M; Vx=\rho_\pi(x)V, \quad x\in N\}.$$
Let $e_N$ be the Jones projection for $E$.
Then using Peter-Weyl theorem we can show
$$\sum_{\pi\in \hG}\sum_{i=1}^{d(\pi)} d(\pi)\Vpi^* e_N \Vpi=1.$$
This means that we can identify $\Xi$ in Theorem 3.3 with $\hG$, and when
$\xi \in \Xi$ and $\pi \in \hG$ are identified we can identify $\Hx$ with
$\cH_\pi$ as well.
Note that $a_\pi=1$ because
$$\eqalign{ (\Vpi, V(\pi)_j)&
=d(\pi)\int_G\alpha_g(\Vpi V(\pi)_j^*)dg
=d(\pi)\sum_{k,l}(\int_G \pi(g)_{k,i}\overline{\pi(g)_{l,j}}dg)
V(\pi)_k V(\pi)_l^* \cr
&=\delta_{i,j}\sum_k V(\pi)_kV(\pi)_k^*=\delta_{i,j}1.}$$

\proclaim Theorem 3.15. Let $G$ be a compact group and $\alpha$ a minimal
action of $G$ on $M$.
Then the map $H \mapsto M^H$ gives one-to-one correspondence between the
lattice of all closed subgroups of $G$ and that of all intermediate
subfactors of $M\supset M^G$. \par

To prove the theorem we need the following lemma, which is
essentially contained in [R2]. For the sake of completeness we give a
proof.

\proclaim Lemma 3.16. Let $G$ be a compact group and $Rep(G)$ the category
of finite dimensional unitary representations of $G$.
For $\pi \in Rep(G)$ $H_\pi$ denotes the representation space of $\pi$.
Suppose we have a Hilbert subspace $K_\pi\subset H_\pi$ for each
$\pi\in Rep(G)$ satisfying the following:
$$K_\pi\oplus K_\sigma \subset K_{\pi \oplus \sigma}, \quad
\pi,\sigma \in Rep(G),$$
$$K_\pi\otimes K_\sigma \subset K_{\pi \otimes \sigma}, \quad
\pi,\sigma \in Rep(G),$$
$$\overline{K_\pi}=K_{\overline{\pi}}, \quad \pi \in Rep(G),$$
where $\overline{\pi}$ is the complex conjugate representation and
$\overline{K_\pi}$ is the image of $K_\pi$ under the natural map from $H_\pi$
to its complex conjugate Hilbert space.
Then there exists a closed subgroup $H \subset G$ such that
$$K_\pi=\{\xi \in H_\pi; \pi(h)\xi=\xi,\quad h\in H\}.$$
\par

\noindent{\it Proof.} Let $B_0$ be the linear span of
$$\{<\pi(\cdot)\xi|\eta>\in C(G); \xi \in K_\pi, \quad \eta \in H_\pi, \quad
\pi \in Rep(G)\}, $$
where $C(G)$ is the C$^*$-algebra of the continuous functions on $G$.
Then by assumption, $B_0$ is a unital *-subalgebra of $C(G)$ that is
globally invariant under the left translation by $G$.
Let $B$ be the norm closure of $B_0$.
Then thanks to [AHKT, Appendix A], there exists a closed subgroup $H\subset G$
such that $B=C(G/H)$.
This implies that $K_\pi$ is a subspace of the set of $H$ invariant vectors
$L_\pi$.
Suppose $\xi \in L_\pi\ominus K_\pi$ and set $f_\eta(g)=<\pi(g)\xi|\eta>$ for
$\eta \in H_\pi$, $g\in G$.
Then $f_\eta\in C(G/H)$.
On the other hand, Peter-Weyl theorem shows that $f_\eta$ is perpendicular to
$C(G/H)$ in $L^2(G)$ because $B_0$ is dense in $C(G/H)$ in uniform norm and
consequently in $L^2(G)$ as well.
Thus $f_\eta=0$ for all $\eta\in H_\pi$ and $\xi=0$.
This proves the statement. Q.E.D. \par

\noindent{\it Proof of Theorem 3.15.} We may assume that $M^G$ is infinite
because after getting the result for $M\otimes B(\ell^2({\bf N}))$ we can
remove $B(\ell^2({\bf N}))$.
It easily follows from the existence of $\{\cH_\pi\}_{\pi \in \widehat{G}}$
that the map is injective.
Let $L$ be an intermediate subfactor and set $\cK_\pi=L\cap \cH_\pi$.
We arrange the orthonormal basis $\{\Vpi\}_{i=1}^{d(\pi)}$ such that
$\{\Vpi\}_{i=1}^{m_\pi}$ is an orthonormal basis of $\cK_\pi$.
Thanks to Lemma 3.8 and Theorem 3.9 $L$ is characterized by
$$\eqalign{L
&=\{x\in M;E(\cK_\pi^\perp x)=0, \quad \pi\in \hG\}\cr
&=\{x\in M; x(\pi)_i=0, \quad i>m_\pi, \quad \pi\in \hG\}. }$$
Thus it is enough to show that there exists a closed subgroup $H\subset G$
such that
$$\cK_\pi=\{V\in \cH_\pi; \alpha_h(V)=V, \quad h\in H\}.$$
Indeed, since $\{\cK_\pi\}_{\pi \in \hG}$ satisfies the assumption of
Lemma 3.8, it is routine to show that one can extend the assignment
$\pi \mapsto \cK_\pi$ to the whole category of representations such that the
assumption of Lemma 3.16 is fulfilled.
Thus Lemma 3.16 captures the desired closed subgroup $H$. Q.E.D.

\noindent{\it Remark 3.17.} It follows from [R1,AHKT] that if H is
a closed subgroup of $G$ as in Theorem 3.15, then $H$ is equal
to the group of all automorphisms of $M$ leaving $M^H$ pointwise
fixed, therefore we have a complete Galois correspondence.
\vskip15pt

\noindent{\bf 4. Kac algebra case.}

In this section we generalize Theorem 3.13 and Theorem 3.15 to the case of
minimal actions of compact Kac algebras.
It turns out that the Galois correspondence holds between the lattice of
intermediate subfactors and that of left coideal von Neumann
subalgebras. We also prove a bicommutant type theorem between the left coideal
von Neumann subalgebras of a compact Kac algebra and right coideal von Neumann
subalgebras of its dual Hopf algebras.

Let $\cA$ be a compact Kac algebra [ES][BS] with coproduct $\delta$, antipode
$\kappa$, and normalized Haar measure $h$, which is a normal trace state.
We regards $\cA$ as a concrete von Neumann algebra represented on the G.N.S.
Hilbert space $\LA$ of $h$ with the G.N.S. cyclic vector $\Oh$.
The multiplicative unitary associated with $\cA$ is defined by
$$V(x \Oh \otimes \xi)=\delta(x)(\Oh\otimes \xi), \quad \xi\in \LA, \;
x\in \cA. \eqno (4.1)$$
Following [BS], we adopt the dual Hopf algebra [BS] rather than the dual Kac
algebra [ES] as the dual object of $\cA$;
the dual Hopf algebra $\dA$ of $\cA$ is the von Neumann algebra generated by
$$\{(id \otimes \omega)(V);\omega \in B(\LA)_*\}$$
with the comultiplication and the antipode,
$$\dd(y)=V^*(1\otimes y)V, \quad \dk(y)=J_{\cA}y^*J_{\cA},\quad
y\in \dA, \eqno (4.2)$$
where $J_{\cA}$ is the canonical conjugation of $\cA$ with respect to $\Oh$.
Let $U\in B(\LA)$ be the unitary operator defined by
$$Ux\Oh=\kappa(x)\Oh, \quad x\in \cA$$
and set
$$\dV=F(U\otimes 1)V(U\otimes 1)F\in \cA\otimes \dA',\eqno (4.3)$$
$$\widetilde{V}=F(1\otimes U)V(1\otimes U)F\in \cA'\otimes \dA,\eqno (4.4)$$
as in [BS] where $F$ is the flip operator of $\LA\otimes \LA$.
$\dV$ and $\widetilde{V}$ are  multiplicative unitaries satisfying
$$\dV^*(\xi \otimes x\Oh)=\delta(x)(\xi \otimes \Oh), \quad
\xi \in \LA,\; x\in \cA, \eqno (3.5)$$
$$\widetilde{V}(y\otimes 1)\widetilde{V}^*=\dd(y), \quad y\in \dA.
\eqno (4.6)$$

\par

A finite dimensional unitary corepresentation $\pi$ is a pair of a finite
dimensional Hilbert space $H_\pi$ and a linear map
$\Gamma_\pi :H_\pi \longrightarrow H_\pi \otimes \cA$ satisfying
$$(\Gamma_\pi \otimes id) \cdot \Gamma_\pi
=(id \otimes \delta)\cdot \Gamma_\pi$$
and the following unitarity condition:
if $\{ e(\pi)_i\}$ is an orthonormal basis of $H_\pi$ and
$$\Gamma_\pi(e(\pi)_j)=\sum_i e(\pi)_i\otimes u(\pi)_{i,j},$$
then $u(\pi)=(u(\pi)_{i,j})$ is unitary as an element in
$M(d(\pi), {\bf C})\otimes \cA$, where $d(\pi)$ is the dimension of $H_\pi$.
We abuse the notation and call $u(\pi)$ a unitary corepresentation as well.
Basic notions such as tensor product, direct sum, complex conjugate
corepresentations, and irreducibility are defined by a standard procedure.
Note that since $\cA$ is a Kac algebra the complex conjugate corepresentation
$u(\overline{\pi})=(u(\overline{\pi})_{i,j}=u(\pi)_{i,j}^*)$ of $u(\pi)$ is
always unitary [W].
Let $\pi$, $\sigma$ be unitary corepresentations of $\cA$.
Then the following orthogonality relation holds:
$$h(u(\pi)_{i,j}^*u(\sigma)_{k,l})={1\over d(\pi)}\delta_{i,k}\delta_{j,l}.$$
Let $\Xi$ be a complete system of representatives of the irreducible
corepresentations of $\cA$.
Then the linear span of $\{u(\pi)_{i,j}\}_{1\leq i,j \leq d(\pi), \;
\pi \in \Xi}$ is a dense in $\cA$ in weak topology.
For $x\in \cA$ we define $x(\pi)_{i,j}$ by
$$x(\pi)_{i,j}=d(\pi)h(u(\pi)_{i,j}^*x).$$
$\{x(\pi)_{i,j}\}$ determines $x$ in the sense that
$x=\sum x(\pi)_{i,j}u(\pi)_{i,j}$ holds in Hilbert space topology in $\LA$.

\proclaim Definition 4.1. A unital von Neumann subalgebra $\cB$ of a Kac
algebra $\cA$ is called a left (right) coideal von Neumann subalgebra
if and only if $\delta(\cB)\subset \cA\otimes \cB$ (respectively
$\delta(\cB)\subset \cB\otimes \cA$) holds. \par

Let $Corep(\cA)$ be the category of finite dimensional unitary
corepresentations of $\cA$.

\proclaim Proposition 4.2. Let $\cA$ be a compact Kac algebra.
Then there exists one-to-one  correspondence between the following two sets.
\item {(i)} The sets of left coideal von Neumann subalgebras of $\cA$.
\item {(ii)} The set of systems of Hilbert subspaces $K_\pi \subset H_\pi$,
$\pi \in Corep(\cA)$ satisfying  the following:
$$K_\pi \oplus K_\sigma \subset K_{\pi \oplus \sigma},
\quad \pi, \sigma \in Corep(\cA).$$
$$K_\pi \otimes K_\sigma \subset K_{\pi \otimes \sigma},
\quad \pi, \sigma \in Corep(\cA).$$
$$\overline{K_\pi}=K_{\overline \pi}, \quad \pi \in Corep(\cA).$$ \par
{\it The correspondence is given as follows.
Let $\{K_\pi\}$ be a system of subspaces satisfying the condition in (ii) and
$\{e(\pi)_i\}_{i=1}^{d(\pi)}$ an orthonormal basis of $H_\pi$ such that
$\{e(\pi)_i\}_{i=1}^{m_\pi}$ is an orthonormal basis of $K_\pi$.
Then the corresponding left coideal von Neumann subalgebra $\cB$ is the
weak closure of the linear span of $\{u(\pi)_{i,j}\}$ $1\leq i \leq d(\pi)$,
$1\leq j\leq m_\pi$, $\pi \in Corep(\cA)$. }\par

\noindent {\it Proof.} First we note that two distinct von Neumann subalgebras
$\cB_1$ and $\cB_2$ give rise to distinct Hilbert subspaces
$\overline{\cB_1\Oh}$, $\overline{\cB_2\Oh}$ because $h$ is a
faithful normal traces.
It is easy to show that the weakly closed linear subspace $\cB$ defined
in the statement is actually a left coideal von Neumann subalgebra,
so it suffices to prove that every left coideal von Neumann subalgebra
$\cB$ arises in this way.
Let $\{e(\pi)_i\}_{i=1}^{d(\pi)}$ be an orthonormal basis of $H_\pi$ and
we set
$$K_\pi=span\{\sum_{j=1}^{d(\pi)}x(\pi)_{i,j}e(\pi)_j; \; x\in \cB,
\quad 1\leq i\leq d(\pi)\}.$$
Since $K_\pi$ does not depend on the choice of the basis, we may and do assume
that $\{e(\pi)_i\}_{i=1}^{m_\pi}$ is an orthonormal basis.
Thus $x(\pi)_{i,j}=0$ for $x\in \cB$, $j>m_\pi$.
We show that $u(\pi)_{i,j} \in \cB$ for $1\leq i\leq d(\pi)$, $1\leq j\leq
m_\pi$.
By the definition of $K_\pi$, for $j$ with $1\leq j\leq m_\pi$ there exist
$x^1,x^2, \cdots x^{d(\pi)} \in \cB$ such that
$\sum_{i=1}^{d(\pi)}x^i(\pi)_{i,k}=\delta_{j,k}$.
Using unitarity of $u(\pi)$ and $\delta(u(\pi)_{p,q})
=\sum_r u(\pi)_{p,r}\otimes u(\pi)_{r,q}$, we get
$$u(\pi)_{i,k}^*\otimes 1=\sum_p(1\otimes u(\pi)_{k,p})
\delta(u(\pi)_{i,p}^*).$$
Since $\cB$ is a left coideal we obtain
$$\eqalign{
\cB\ni &\sum_i(h\otimes id)((u(\pi)_{i,k}^*\otimes 1)\delta(x^i))
=\sum_{i,p}(h\otimes id)((1\otimes u(\pi)_{k,p})\delta(u(\pi)_{i,p}^*x^i))\cr
&=\sum_{i,p}x^i(\pi)_{i,p}u(\pi)_{k,p}=u(\pi)_{k,j}.}$$
Thus $\cB$ is characterized as
$$\cB=\{x\in \cA;\; x(\pi)_{i,j}=0, \quad \pi \in \Xi, \; j>m_\pi\}.$$
Since $\cB$ is a *-subalgebra, the natural extension of
$\{K_\pi\}_{\pi \in \Xi}$ to the whole category of unitary corepresentations
satisfies the three conditions of (ii). Q.E.D.

\proclaim Definition 4.3. Let $\Gamma : M\longrightarrow M\otimes \cA$ be an
action of a compact Kac algebra $\cA$ on a factor $M$.
\item {(i)} $\Gamma$ is called minimal if and only if the linear span of
$\{(\omega \otimes id)\cdot \Gamma (M); \omega \in M_*\}$ is dense in $\cA$
and the relative commutant of the fixed point algebra
$M^\Gamma=\{x\in M; \Gamma(x)=x\otimes 1\}$ in $M$ is trivial.
\item {(ii)} Let $\cB$ be a left coideal von Neumann subalgebra of $\cA$.
The intermediate subalgebra $M(\cB)$ of $M^\Gamma \subset M$ associated
to $\cB$ is defined by
$$M(\cB)=\{x\in M; \Gamma(x)\in M\otimes \cB\}.$$

\proclaim Theorem 4.4. Let $\Gamma :M\longrightarrow M\otimes \cA$ be
a minimal action of a compact Kac algebra $\cA$ on a factor $M$.
Then the map $\cB \mapsto M(\cB)$ gives one-to-one correspondence
between the lattice of left coideal von Neumann subalgebras of $\cA$
and that of  the intermediate subfactors of $M^\Gamma \subset M$.
\par

\noindent{\it Proof.} For the same reason as in the proof of Theorem 3.15
we may assume that $M^\Gamma$ is infinite.
Note that there exists a normal conditional expectation $E\in \ce M{M^\Gamma}$
given by
$$E(x)\otimes 1=(id\otimes h)\cdot \delta(x), \quad x\in M.$$
In exactly the same way as in the case of compact group actions,
for each $\pi \in \Xi$ one can find a Hilbert space $\cH_\pi$ in $M$
with support 1 and its basis $\{V(\pi)\}_{i=1}^{d(\pi)}$ satisfying
$$\delta(V(\pi)_i)=\sum_jV(\pi)_j\otimes u(\pi)_{j,i}.$$
Thus, as before we can identify our $\Xi$ with that in Theorem 3.3 and
we get $a_\pi=1$ thanks to the orthogonality relation.
Let $L$ be an intermediate subfactor and $\cK_\pi=L\cap \cH_\pi$.
Thanks to Lemma 3.8 and Theorem 3.9, $L$ is generated by $M^\Gamma$ and
$\{\cK_\pi\}_{\pi \in \Gamma}$, and is characterized by
$$L=\{x\in M; \; E(\cK_\pi^\perp x)=0, \quad \pi \in \Xi\}. $$
Therefore, as in the proof of Theorem 3.16 we can conclude $L=M(\cB)$
by using Proposition 4.2, where $\cB$ is the left coideal von Neumann
subalgebra corresponding to $\{\cK_\pi\}_\pi$.
The map is injective because 2 distinct systems of subspaces satisfying the
assumption of Proposition 4.2 (ii) give rise to 2 distinct intermediate
subfactors. Q.E.D.

\noindent{\it Remark 4.5}. The crossed product $M\times_\Gamma \dA$ is the
von Neumann algebra generated by $\Gamma(M)$ and ${\bf C}\otimes \dA$.
As is expected, we can identify the basic extension $M_1$ with
$M\times_\Gamma \dA$ if the action is minimal as follows.
Let $e_0$ be the projection in $\dA$ corresponding to the trivial
corepresentation of $\cA$ and we set $e=1\otimes e_0$.
Since we have the dual operator valued weight of the crossed product whose
restriction to $\dA$ is a semi-finite trace (Plancherel weight),
if $MeM$ is dense in $M\times_\Gamma \dA$ we can apply Lemma 2.4 and get the
result.
Indeed, it is known [BS] that $\dA$ is a direct sum of type I$_{d(\pi)}$
factors $\dA_\pi$, $\pi \in \Xi$ and the multiplicative unitary $V$ can be
expanded as
$$V=\sum_{\pi\in \Xi}\sum_{1\leq i,j \leq d(\pi)} e(\pi)_{i,j}
\otimes u(\pi)_{i,j}, $$
where $\{e(\pi)_{i,j}\}$ are matrix units  of $\dA_\pi$.
Thanks to (4.1), we have
$$e(\pi)_{i,j}u(\sigma)_{k,l}\Oh
=\delta_{\pi, \sigma}\delta_{j, l}u(\pi)_{k,i}\Oh.$$
Now, we show
$$d(\pi) \delta(V(\pi)_i^*)e\delta(V(\pi)_j)
=1\otimes \widehat{\kappa}(e(\pi)_{j,i}).$$
From $\delta(V(\pi)_i)=\sum_k V(\pi)_k\otimes u(\pi)_{k,i}$ We get
$$\delta(V(\pi)_i^*)e\delta(V(\pi)_j)
=\sum_k 1\otimes u(\pi)_{k,i}^*e_0u(\pi)_{k,j}.$$
Thanks to the orthogonality relation, we obtain
$$d(\pi) \sum_k u(\pi)_{k,i}^*e_0u(\pi)_{k,j}u(\sigma)_{p,q}^*\Oh
=\delta_{\pi, \sigma}\delta_{j,q}u(\pi)_{p,i}^*,$$
where we use the fact that $e_0$ is the projection onto the space spanned by
$\Oh$.
On the other hand,
$$\eqalign{\widehat{\kappa}(e(\pi)_{j,i})u(\sigma)_{p,q}^*\Oh
&=J_{\cA}e(\pi)_{i,j}J_{\cA}u(\sigma)_{p,q}^*\Oh
=J_{\cA}e(\pi)_{i,j}u(\sigma)_{p,q}\Oh\cr
&=\delta_{\pi,\sigma}\delta_{j,q}J_{\cA}u(\pi)_{p,i}\Oh
=\delta_{\pi,\sigma}\delta_{j,q}u(\pi)_{p,i}^*\Oh .}$$
Thus $\delta(M)e\delta(M)$ is dense in $M\times_\Gamma \dA$.

The above theorem suggests that it is worth while studying the structure of
the lattice of the left coideal von Neumann subalgebras of Kac algebras.
For compact and discrete Kac algebras we have the following:

\proclaim Theorem 4.6. Let $\cA$ be a compact Kac algebra and $\dA$ its dual
Hopf algebra represented on $\LA$.
Let $\cB \subset \cA$ be a left coideal von Neumann subalgebra and
$\cC\in \dA$ a right coideal von Neumann subalgebra.
Then the following hold:
\item {(i)} $\cB'\cap \dA$ is a right coideal von Neumann subalgebra of
$\dA$ and $(\cB'\cap \dA)'\cap \cA=\cB$.
\item {(ii)} $\cC'\cap \cA$ is a left coideal von Neumann subalgebra of
$\cA$ and $(\cC'\cap \cA)'\cap \dA=\cC$.
\item {(iii)} Set $\tB=\dk(\cB'\cap \dA)$.
Then the map given by $\cB\mapsto \tB$ is a lattice anti-isomorphism between
the set of left coideal von Neumann subalgebras of $\cA$ and that of $\dA$.
\par
\noindent{\it Proof.} (i): Let $E_\cB$ be the $h$ preserving conditional
expectation in $\ce {\cA}{\cB}$ and $e_\cB$ its Jones projection, i.e.
$e_\cB$ is the projection defined by $e_\cB x\Oh=E_\cB(x)\Oh$, $x\in \cA$.
Note that $e_\cB\in \cB'$ and $\{e_\cB\}' \cap \cA=\cB$ hold.
Thus to prove $(\cB'\cap \dA)'\cap \cA=\cB$ it suffices to show
$e_\cB\in \cB'\cap \dA$.
First, we prove $(id \otimes E_\cB)\cdot \delta =\delta \cdot E_\cB$.
Let $\{\cK_\pi\}_{\pi\in \Xi}$ be the system of Hilbert subspaces corresponding
to $\cB$ and $\{e(\pi)_i\}_{i=1}^{d(\pi)}$ an orthonormal basis of $\cH_\pi$
such that $\{e(\pi)_i\}_{i=1}^{m_\pi}$ is a basis of $\cK_\pi$.
As we saw in the proof of Proposition 4.2, the linear span of
$\{u(\pi)_{i,j}\}$,
$\pi \in \Xi$, $1\leq i\leq d(\pi)$, $1\leq j\leq m_\pi$ is dense in $\cB$ in
weak topology.
Let $x\in \cA$, $1\leq j\leq m_\pi$.
Then we get
$$\eqalign {(id \otimes h)((1\otimes u(\pi)_{i,j})\delta(x))
&=\sum_{k=1}^{d(\pi)}u(\pi)_{k,i}^*(id \otimes h)(\delta(u(\pi)_{k,j}x)) \cr
&=\sum_{k=1}^{d(\pi)}u(\pi)_{k,i}^*h(u(\pi)_{k,j}x) \cr
&=\sum_{k=1}^{d(\pi)}u(\pi)_{k,i}^*h(u(\pi)_{k,j}E_\cB(x))\cr
&=(id \otimes h)((1\otimes u(\pi)_{i,j})\delta(E_\cB(x))), }$$
which implies $(id \otimes E_\cB)\cdot \delta (x)=\delta \cdot E_\cB(x)$.
Let $\widehat{V}$ be the multiplicative unitary defined in (4.3).
Then thanks to (4.5), for $\xi \in \LA$ and $x\in \cA$ we get
$$\eqalign{
\widehat{V}^*(1\otimes e_\cB)(\xi \otimes x\Oh)
&=\widehat{V}^*(\xi \otimes E_\cB(x)\Oh)
=\delta (E_\cB(x))(\xi \otimes \Oh) \cr
&=(id \otimes E_\cB)\cdot \delta(x)(\xi \otimes \Oh)
=(1\otimes e_\cB)\widehat{V}^*(\xi\otimes x\Oh),}  $$
and so $(1\otimes e_\cB)$ commutes with $\widehat{V}$.
Since $\{(\omega\otimes id)(\widehat{V});\; \omega \in B(\LA)_*\}$ is dense in
$\dA'$, $e_\cB\in \dA$.
Let $x\in \cB$, $y\in \cB'\cap \dA$. Then
$$\eqalign{
\hat{\delta}(y)(x\otimes 1)
&=V^*(1\otimes y)V(x\otimes 1)
=V^*(1\otimes y)\delta(x)V=V^*\delta(x)(1\otimes y)V \cr
&=(x\otimes 1)V^*(1\otimes y)V=(x\otimes 1)\hat{\delta}(y).}$$
Thus $\cB'\cap \dA$ is a right coideal von Neumann subalgebra of $\dA$.
(ii): In a similar way as above one can show that $\cC'\cap \cA$ is a left
coideal von Neumann subalgebra of $\cA$.
Let $\cC_0=(\cC'\cap \cA)'\cap \dA$.
Then it is easy to show $\cC_0'\cap \cA=\cC'\cap \cA$.
Thus to prove $\cC_0=\cC$, it suffices to prove that if $\cC_1$ and $\cC_2$
are distinct right coideal von Neumann subalgebras of $\dA$, then
$\cC_1'\cap \cA$ and $\cC_2' \cap \cA$ are distinct.
Since the Plancherel weight of $\dA$ is the restriction of the usual trace
on $B(\LA)$, there exists a trace preserving conditional expectation
$F\in \ce {B(\LA)}{\dA}$.
Note that one can identify $F$ with the dual weight of the crossed product of
$\dA$ and $\hat{\hat{\cA}}=\cA'$ when $\hat{\delta}$ is regarded as an action
of $\dA$ on itself.
Thus the restriction of $F$ to $\cA'$ is a trace.
We claim that $F((\cC' \cap \cA)')=\cC$ for every right coideal von Neumann
subalgebra $\cC\subset \dA$.
To prove the claim it is enough to show that $\cC \cdot \cA'$ is weakly
dense in $(\cC \cup \cA')''$ because of $(\cC' \cap \cA)'=(\cC \cup \cA')''$.
Let $\widetilde{V}$ be as in (4.4).
Thanks to (4.4) and (4.6), for $c\in \cC$ and $\omega \in B(\LA)_*$ we get
$$(id\otimes \omega)(\widetilde{V})c
=(id\otimes \omega)(\widetilde{V}(c\otimes 1))
=(id\otimes \omega)(\hat{\delta}(c)\widetilde{V})
\in \overline{\cC \cdot \cA'}^w,$$
which shows $\overline{\cC \cdot \cA'}^w=(\cC \cup \cA')''$.
Using the claim, now we can show that if $\cC_1\neq \cC_2$ are right coideal
von Neumann subalgebras of $\dA$, $(\cC_1' \cap \cA)'\neq (\cC_2' \cap \cA)'$,
and so $(\cC' \cap \cA)' \cap \dA=\cA$. (iii): This is a direct consequence
of (i) and (ii). Q.E.D.

In what follows, we assume $n:=dim \cA<\infty$.
Let $\epsilon$ and $\hat{\epsilon}$ be the counit of $\cA$ and $\dA$, and
$e$ and $\hat{e}$ the integrals of $\cA$ and $\dA$; $e$ and $\hat{e}$ are
the minimal central projections satisfying $ex=e\epsilon(x)$, $x\in \cA$, and
$\hat{e}y=\hat{e}\hat{\epsilon}(y)$, $y\in \dA$.
It is known that the G.N.S. cyclic vector $\Omega_{\hat{h}}$ of the
normalized Haar measure $\hat{h}$ of $\dA$ can be identified with
$\sqrt{n}e\Oh$ and we have $\sqrt{n}\hat{e}\Omega_{\hat{h}}=\Omega_h$ as
well [KP].
The dual pairing between $\cA$ and $\dA$ can be written in terms of the
Hilbert space inner product as follows:
$$<x,y>=\sqrt{n}<x\Oh|y^*\Omega_{\hat{h}}>, \quad x\in \cA\; y\in \dA.
\eqno (4.7) $$

The following is a space free description of the anti-isomorphism of the two
lattices.

\proclaim Proposition 4.8. Let $\cA$ be a finite dimensional Kac algebra and
$\cB$ a left coideal von Neumann subalgebra of $\cA$.
We set
$$\tB=\{y\in \dA;<xb,y>=\epsilon(b)<x,y>, \quad x\in \cA, \; b\in \cB\}.$$
Then the following hold:
\item {(i)} $\tB$ is a left coideal von Neumann subalgebra of $\dA$
with $dim\cB \cdot dim\tB=dim\cA$.
\item {(ii)} $\widetilde{\tB}=\cB$.
\item {(iii)} $\tB=\dk(\cB'\cap \dA)$. \par

\noindent{\it Proof.}
(i): It is routine to show that $\tB$ is a left coideal von Neumann subalgebra
of $\dA$.
Using (4.7) for $x\in \cA$, $b\in \cB$, and $y\in \dA$ we get the following:
$$\eqalign{<xb,y>
&=\sqrt{n}<xb\Oh|y^*\Odh>=\sqrt{n} <J_\cA b^* x^*\Oh|y^*\Odh> \cr
&=\sqrt{n}<b J_\cA y^*\Odh|x^*\Oh>
=\sqrt{n}<b \dk(y)\Odh|x^*\Oh>.}$$
On the other hand we have
$$\epsilon(b)<x,y>=\epsilon(b)\sqrt{n}<x\Oh|y^*\Odh>
=\epsilon(b)\sqrt{n}<\dk(y)\Odh|x^*\Oh>, $$
and so $\tB$ is characterized as
$$\tB=\{y\in \dA; b\dk(y)\Odh=\epsilon(b)\dk(y)\Odh, \quad b\in \cB\}.$$
Let $E_\cB$ and $E_\tB$ be the $h$ and $\dh$ preserving conditional
expectations onto $\cB$ and $\tB$ respectively, and $e_\cB$ and $e_\tB$ the
corresponding Jones projections.
The above characterization shows
$$\tB\supset \dk(\cB'\cap \dA)\ni J_\cA e_\cB J_\cA=e_\cB.$$
More specifically we show $e_\cB=n\epsilon \cdot E_\cB(e)E_\tB(\hat{e})$.
Indeed, using $\dk(e_\cB)=J_\cA e_\cB J_\cA=e_\cB$ we get the following for
$\tb\in \tB$:
$$\eqalign{\dh(e_\cB\tb)
&=\dh(\dk(\tb)e_\cB)=<\dk(\tb)e_\cB\Odh|\Odh> \cr
&=\sqrt{n}<e_\cB e\Oh|\dk(\tb^*)\Odh>=\sqrt{n}<E_\cB(e)\Oh|\dk(\tb^*)\Odh>\cr
&=\sqrt{n}<\Oh|E_\cB(e)\dk(\tb^*)\Odh>
=\epsilon\cdot E_\cB(e)\sqrt{n}<\Oh|\dk(\tb^*)\Odh>\cr
&=\epsilon\cdot E_\cB(e)\sqrt{n}<\dk(\tb)\Oh|\Odh>
=\epsilon\cdot E_\cB(e)\hat{\epsilon}(\dk(\tb))\sqrt{n}<\Oh|\Odh>\cr
&=\epsilon\cdot E_\cB(e)\hat{\epsilon}(\dk(\tb))
=n\dh(\hat{e}\tb)\epsilon\cdot E_\cB(e).} $$
Thus we obtain the claim.
Note that $\dh$ is the restriction of the normalized trace of $B(\LA)$, and
so $\dh(e_\cB)={dim \cB\over n}$.
Thus we get
$$\epsilon \cdot E_\cB(e)={dim \cB\over n}, $$
$$e_\cB=dim \cB E_\tB(\hat{e}).$$
In the same way we can get
$$\hat{\epsilon} \cdot E_\tB(\hat{e})={dim \tB\over n}.  $$
Since $e_\cB$ is a projection,
$$\eqalign{e_\cB=e_\cB^2
&=(dim \cB)^2E_\tB(\hat{e}E_\tB(\hat{e}))
=(dim \cB)^2 \hat{\epsilon}\cdot E_\tB(\hat{e})E_\tB(\hat{e})\cr
&=dim \cB \hat{\epsilon}\cdot E_\tB(\hat{e}) e_\cB .}$$
Therefore, $dim \cB dim \tB=n$.
(ii): It is easy to show $\cB \subset \widetilde{\tB}$.
Since $\tB$ is a left coideal von Neumann subalgebra of $\dA$ we also have
$dim \tB \cdot dim \widetilde{\tB}=n$, and so $\cB =\widetilde{\tB}$.
(iii): In a similar way as in (i),  one can show $dim \cB \cdot
dim (\cB'\cap \dA)=n$.
Since we have the inclusion $\tB\supset \dk(\cB'\cap \dA)$ we get the equality.
Q.E.D.

\noindent{\it Remark 4.9.} Let $\Gamma :M\longrightarrow M\otimes \cA$ be
a minimal action on a factor $M$ and $L$ be an intermediate subfactor of
$M^\Gamma \subset M$.
Thanks to Theorem 4.6 there exists a left coideal von Neumann subalgebra
$\cB \subset \cA$ such that $L=M(\cB)$.
Let $L_1=J_M L' J_M$, which is an intermediate subfactor of $M\subset M_1$.
Under the identification of $M_1$ with $M\times_\Gamma \dA$, we actually have
$L_1=M_1(\tB)$.
Let $e_L$ be the Jones projection for $L$.
Then,
$$\eqalign{L_1
&=(M\cup \{e_L\})''=(M\cup (L_1\cap N'))''\cr
&=(M\cup J_M(L'\cap M_1)J_M)''=(M\cup j(L'\cap M_1))'',}$$
where $j$ is the anti-automorphism of $N'\cap M_1$ defined by
$j(x)=J_M x^*J_M$, $x\in N'\cap M$.
It is known [D] that under our identification, $N'\cap M_1$ is identified
with ${\bf C}\otimes \dA$ and $j$ is identified with $\dk$, so
$j(L' \cap M_1)$ is identified with ${\bf C}\otimes \tB$.
Thus $L_1$ is identified with the intermediate subfactor generated by
$\Gamma(M)$ and ${\bf C}\otimes \tB$, which proves the remark.

\vskip20pt
\centerline{\bf References}
\medskip

\item { }[AHKT] H. Araki, R. Haag, D. Kastler, M. Takesaki, {\it Extension of
KMS ststes and chemical potential,} Commun. Math. Phys. {\bf 53} (1977),
97-134.
\item { }[B] M. W. Binder, {\it Induced factor representations of discrete
groups and their types,} J. Funct. Anal. {\bf 115} (1993), 294-312.
\item { }[BS] S. Baaj, G. Skandalis, {\it Unitaires multiplicatifs et
dualit{\'e} pour les produits crois{\'e}s de $C^*$-algebres,} Ann. Scient.
{\`E}c. Norm. Sup. {\bf 26} (1993), 425-488.
\item { }[Ch] H. Choda. {\it A Galois correspondence in a von Neumann
algebra}, Tohoku Math. J. {\bf 30} (1978), 491-504.
\item { }[C] A. Connes, {\it On spatial theory of von Neumann algebras,}
J. Funct. Anal. {\bf 35} (1980), 153-164.
\item { }[C2] A. Connes, {\it Non-commutative Geometry},
Academic Press (1994).
\item { } [CT] A. Connes, M. Takesaki, {\it The flow of weights of on
factors of type III,} Tohoku Math. J. {\bf 29} (1977), 473-575.
\item { }[D] M.-C. David, {\it Paragroup d'Adrian Ocneanu et algebra de Kac,}
Pac. J. Math. (to appear).
\item { }[DR] S. Doplicher, J.E. Roberts {\it Why there is a field
algebra with a compact gauge group describing the superselection
strucure in particle physics} Commun. Math. Phys. {\bf 131} 51-107
(1990). \item { }[EN] M. Enock, R. Nest, {\it Irreducible inclusions
of factors  and multiplicative unitaries,} preprint.
\item { }[ES] M. Enock, J.-M. Schwartz, {\it Kac algebras and duality of
locally compact groups,} Springer, Berlin, 1992.
\item { }[GHJ] F. Goodman, P. de la Harpe, V. Jones, {\it Coxeter
graphs and towers of algebras,} MSRI Publications 14, Springer Verlag,
1989.
\item { }[H1] U. Haagerup, {\it Operator valued weights in von Neumann
algebras, I,} J. Funct. Anal. {\bf 32} (1979), 175-206.
\item { }[H1] U. Haagerup, {\it Operator valued weights in von Neumann
algebras, II,} J. Funct. Anal. {\bf 33} (1979), 339-361.
\item {  }[HK] T. Hamachi, H. Kosaki, {\it Orbital factor map,}
Erg. Th. Dyanam. Syst.
\item { }[HS] U. Haagerup, E. St{\o}rmer, {\it Subfactors of a factor of type
III$_\lambda$ which contain a maximal centralizer,} Internat. J. Math. {\bf 6}
(1995), 273-277.
\item { }[Hi] F. Hiai, {\it Minimizing indices of conditional
expectations on a subfactor,} Publ. RIMS, Kyoto Univ. {\bf24} (1988), 673-678.
\item { }[HO] R. Hermann, A. Ocneanu, {\it Index theory and Galois theory for
infinite index inclusions of factors,} C.R. Acad. Sci. Paris, {\bf 309}
(1989), 923-927.
\item { }[I1] M. Izumi, {\it Application of fusion rules to classification of
subfactors,} Publ. RIMS, Kyoto Univ. {\bf 27} (1991),953-994.
\item { }[I2] M. Izumi, {\it Subalgebras of infinite C$^*$-algebras with
finite Watatani indices. II. Cuntz-Krieger algebras,} preprint.
\item { }[J] V. Jones, {\it Index for subfactors,} Invent. Math. {\bf
72} (1983), 1-25.
\item { }[K] A. Kishimoto, {\it Remarks on compact automorphism groups of a
certain von Neumann akgebra,} Publ. RIMS, Kyoto Univ. {\bf 13} (1977), 573-581.
\item { }[Ko1] H. Kosaki, {\it Extension of Jones theory on index to
arbitrary factors,} J. Funct. Anal. {\bf 66} (1986), 123-140.
\item { }[Ko2] H. Kosaki, {\it Characterization of crossed product
(properly infinite case),} Pacif. J. Math. {\bf 137} (1989),
159-167.  \item { }[KL] H. Kosaki, R. Longo, {\it A remark on the
minimal index  of subfactors,} J. Funct. Anal. {\bf 107} (1992),
458-470.  \item { }[KY] H. Kosaki, S. Yamagami, {\it Irreducible
bimodules associated  with crossed product algebras,} Internat. J.
Math. {\bf 3} (1992), 661-676.  \item { }[KP] G.I. Kac, V.G.
Pljutkin, {\it Finite ring groups,} Trans.  Moscow Math. Soc.
(1966), 251-294, Translated from Trudy Moskov. Mat. Obsc.  {\bf 15}
(1966), 224-261.  \item { }[L1] R. Longo, {\it Index of subfactors
and statistics of  quantum fields I ,} Commun. Math. Phys. {\bf 126}
(1989), 217-247. \item { }[L2] R. Longo, {\it Index of subfactors
and statistics of  quantum fields II ,} Commun. Math. Phys. {\bf
130} (1990), 285-309. \item { }[L3] R. Longo, {\it Simple injective
subfactors}, Adv. Math.  {\bf 63} (1987), 152-171.
\item { }[L4] R. Longo, {\it Minimal index and braided subfactors,} J.
Funct. Anal.
 {\bf 109} (1992), 98-112.
\item { }[N] Y. Nakagami, {\it Essential spectrum $\Gamma(\beta)$ of
a dual action on a von Neumann algebra,} Pacif. J. Math. {\bf 70} (1977),
437-478.
\item { }[NT] N. Nakamura, Z. Takeda, {\it A Galois theory for finite
factors,} Proc. Japan Acad. {\bf 36} (1960), 258-260.
\item { }[NTs] Y. Nakagami, M. Takesaki {\it ``Duality for crossed product
of von Neumann algebras"} Lecture Notes in Math. {\bf 731}, 1979,
Springer-Verlag, Berlin-Heidelberg-New York.
\item { }[O] A. Ocneanu, {\it Quantized group string algebra and
Galois  theory for algebra,} in ``Operator algebras and
applications, Vol. 2  (Warwick, 1987),'' London Math. Soc. Lect.
Note Series Vol. 136,  Cambridge University Press, 1988, pp. 119-172.
\item { }[P] S. Popa, {Classification of amenable subfactors of type II,}
Acta Math. {\bf 172} (1994), 352-445.
\item { }[PP1] M. Pimsner, S. Popa, {\it Entropy and index for
subfactors,}
Ann. sient. $\acute {\rm E}$c. Norm. Sup. {\bf 4}  57- 106, 1986.
\item { }[PP2] M. Pimsner, S. Popa, {\it Iterating basic construction,}
Trans. Amer. Math. Soc. {\bf 210} (1988), 127-133.
\item { }[S] {\c S}. Str{\v a}til{\v a}, {\it Modular theory in operator
algebras, } Editra Academiei and Abacus Press 1981.
\item { }[R1] J. E. Roberts, {\it Cross product of von Neumann
algebras by group duals,} Symposia Math. Vol. XX (1976).
\item { }[R2] J. E. Roberts, {\it Spontaneously broken gauge
symmetries and superselection rules,} Proc. of ``School of
Mathematical Physics", Universita' di Camerino, 1974.
\item { } [Su] N. Suzuki {\it Crossed product of rings of operators}
T\^ohoku Math. J. {\bf 11} 113-124 (1960).
\item { }[Sy] W. Szyma\'nski, {\it Finite index
subfactors and Hopf algebra  crossed products,} Proc. Amer. Math.
Soc., {\bf 120} (1994), 519-528.  \item { }[W] S. L. Woronowicz,
{\it Compact matrix pseudogroups,} Commun.  Math. Phys., {\bf 111}
(1987), 618-665.  \item { }[Y1] S. Yamagami, {\it Modular theory for
bimodules,} J. Funct. Anal.  {\bf 125} (1994), 327-357.
\item { }[Y2] S. Yamagami, {\it On Ocneanu's characterization of crossed
products,} preprint.

\bigskip

\noindent {\bf Acknowledgments.} We wish to thank J.E.
Roberts for conversations. M.I. and S.P. would like to thank
the CNR and the Universit\`a di Roma ``Tor Vergata'' for supporting
several visits to Rome that made this collaboration possible.
\bigskip\noindent
\input amstex
\documentstyle{amsppt}
\NoBlackBoxes
\magnification 1200
\def\Aut{\hbox{ Aut }}
\def\dim{\hbox{ dim }}
{\bf Appendix} ({\it Added October 30, 1997})\par
\medskip
Related to Theorem 3.3 and Remark 3.4, we give here an example of an
irreducible inclusion of factors $N\subset M$ with a normal
conditional expectation $E$ such that $N\subset M$ is discrete,
$E^{-1}$ is a semifinite trace on $N'\cap M_1$ (so that
$B_1=B_2=C=\{0\}$) yet $E^{-1}\circ j \neq E^{-1}$. In fact,
our factors $N, M$ are hyperfinite of type II$_1$
with $E\in{\Cal E}(M,N)$ being  the unique normal conditional expectation
preserving  the trace
$\tau$ on $M$ and $E^{-1}$ being a semifinite trace on $N'\cap M_1$. Thus,
while the irreducibility of an inclusion of (type II$_1$) factors $N\subset M$
with $[M:N]<\infty$ automatically entails its extremality (thus, the
trace-preservingness of $j=J_M\cdot J_M$
on $N'\cap M_1)$ this is no
longer the case when $[M:N]=\infty,$ even if $N\subset M$ is discrete.

Our construction is based on Powers
binary shifts and their properties ([Po], [PoPr]).
\bigskip
\proclaim {Lemma A.1 ([PoPr])} Let $\sigma$ be a bilateral Powers binary
shift acting
on $\{u_n\}_{n\in\Bbb Z}$ as in [Po], such that each half-line bitstream of
$\sigma$ is aperiodic.\par
Let $P= vN\{u_n\}_{n\in\Bbb Z}$ and $N= vN\{u_n\}_{n\ge 0}.$ Then the following
hold true:

\item{(i)} $N$ and $P$ are factors;\par
\item{(ii)} $\sigma(N)\subset N$ and $[N:\sigma(N)]<\infty.$\par
\item{(iii)} $\sigma^n(N)'\cap N=\Bbb C, \forall n\ge 1.$\par
\item{(iv)} $\bigcap\limits_{n\ge 1}^{}\sigma^n(N)=\Bbb C 1.$\par
\item{(v)} $\bigcup\limits^{}_{n\ge 1}\sigma^n(N)$ is a dense $^*$-subalgebra
of $P.$\endproclaim
\demo{Proof} All these are well known properties from [Po] [PoPr].\enddemo
\hfill{Q.E.D.}
\bigskip
\proclaim{ Proposition A.2} Let $P$ be a type II$_1$ factor with an
aperiodic automorphism
$\sigma\in \Aut P$ and a subfactor $N\subset P$ such that $P, \sigma, N$
satisfy the conditions (i)-(v) of the previous Lemma.\par
Let $M=P\rtimes_\sigma\Bbb Z.$ Then we have:

a) $N'\cap M=\Bbb C1.$\par
b) $N \subset M$ is discrete, i.e.,
$L^2( M_1 ,  {\text{\rm Tr}})$ is generated by $N-N$ sub-bimodules which
have finite dimension
both as left and right $N-$modules, where {\rm Tr} is the unique
semifinite trace on $M_1=\langle N,M\rangle$ such that {\rm Tr}$e_N
=1$.\par
c) $J_M\cdot J_M$ is not {\rm Tr}-preserving on $N'\cap M_1 ,$
equivalently, there are
irreducible $N-N$ sub-bimodules of $L^2( M_1 ,  {\text{\rm Tr}})$ for which
the left dimension
over $N$ does not coincide with the right dimension over $N.$
\endproclaim
\demo{Proof} a). By property (iii) we have $N'\cap\sigma^{-n}(N)=\Bbb C1.$
Thus, if $a\in N'\cap P$ then $\|E_{\sigma^{-n}(N)} (a)-a\|_2\to 0$ (by (ii)
and (v)) and $E_{\sigma^{-n}(N)} (a)\in N'\cap \sigma^{-n}(N)=\Bbb C$ (by
commuting squares). Thus $a\in\Bbb C1,$ showing that $N'\cap P=\Bbb C.$
Similarly $\sigma^n(N)'\cap P=\Bbb C,\, \forall n\ge 1.$\par
Assume now that $a=\sum\limits^{}_{n\in\Bbb Z}b_n u^n$ satisfies $ax=xa,
\,\forall\, x\in N.$ Thus, if $b_n\ne 0$ for some $n$ then $xb_n=b_n\sigma^n(x)
\forall x\in N.$ By using that $\sigma^n(N)\subset N,$ it follows that
$xb_nb_n^*=b_n\sigma^n(x)b^*_n=b_nb^*_nx,\, \forall\, x\in N.$ Thus $b_nb^*_n\in
\Bbb C1$ so that $b_n$ is a (multiple of a) unitary element $v\in P$
satisfying
$$
x v=v\sigma^n(x),\, \forall\, x\in N.\eqno (1)$$
In particular, we have
$$
\sigma^n(x)v=v\sigma^{2n}(x), \,\forall\, x\in N\eqno (2)
$$
Also, by applying $\sigma^n$ to both sides of (1) we get
$$
\sigma^n(x)\sigma^n(v)=\sigma^n(v)\sigma^{2n}(x),\, \forall\, x\in
N.\eqno (3)
$$From (2) and (3) we get:
$$
v^*\sigma^n(x)v=\sigma^n(v^*)\sigma^n(x)\sigma^n(v), \,\forall\, x\in
N.\eqno (4)
$$
Thus, $\sigma^n(v)=\alpha v$ for some $\alpha\in\Bbb C1.$ Let then $m_k
\nearrow\infty$ be such that $\alpha^{m_k}\to 1.$ Then $\|\sigma^{nm_k}(v)
-v\|\to 0$ as $k\to\infty.$ But $\bigcap\limits^{}_{m\ge 0}\sigma (N)=
\Bbb C1$ clearly implies $\sigma$ is mixing, i.e. $\underset {u\to\infty}
\to\lim\tau(\sigma^n(x)y)=\tau(x)\tau (y), \,\forall\, x,y\in P.$ A
contradiction unless $v\in\Bbb C1,$ showing that $N'\cap M=\Bbb C.$\par
\medskip
(b). Let $K_{n,m}=u^{-n}L^2(N)u^{n+m},$ for $n\ge 0, m\in\Bbb Z, n\ge m.$
It is trivial to see that $K_{n,m}\nearrow L^2(P)u^m,$ as $n\nearrow\infty.$
Thus $\vee\{ K_{n,m}| n\geq m,n\geq 0,m\in\Bbb Z\} =L^2(M),$ with all
$K_{n,m}$ being $N-N$ bimodules.\par
Also, since as a left $N$ module $K_{n,m}=L^2(\sigma^{-n}(N))u^m$ is isomorphic
to $L^2(\sigma^{-n}(N)),$ we have $\dim({}_N K_{n,m})=[N:\sigma(N)]^n<\infty.$
Furthermore, as a right $N$-module $K_{n,m}=u^mL^2(\sigma^{-n-m}(N))$ is
isomorphic to $L^2(\sigma^{-n-m}(N)),$ so that we have $\dim
(K_{n,m N})=
[N:\sigma(N)]^{n+m}<\infty.$\par
This shows that $N\subset M$ is discrete.\par
\medskip
c). This part is now clear, since we showed above that there exist
sub-bimodules
$K\subset L^2( M_1 ,  {\text{\rm Tr}})$ which are finitely generated both
as left and right
modules, but with different corresponding dimensions, (e.g., just take
$K=K_{n,m}$ for some $n\ge m, n\ge 0, m\ne 0.$)\par
\hfill{Q.E.D.}
\bigskip
\proclaim{Corollary A.3}
There exist irreducible discrete inclusions of hyperfinite type II$_1$ factors
$N\subset M$ for which $J_M\cdot J_M$ is not trace preserving on
$N'\cap  M_1 ,$ equivalently for which
$ {\text{\rm Tr}}_{ M_1 }$ and $ {\text{\rm Tr}}_{N'}$ do not agree on
$N'\cap M_1$
or further, for which the local indices $\ [pM_1 p:Np]$ are not
equal to $({\text{\rm Tr}}\, p)^2$ for all $p\in N'\cap  M_1$.
\endproclaim
\medskip
{\bf References}
\medskip\noindent
[Po] R. T. Powers
{\it An index theory for semigroups of *-endomorphisms of $B(H)$ and
type II$_1$ factors},  Can. J. Math {\bf 40} (1988), 86-114.
\par\noindent
[PoPr] R. T. Powers, G. L. Price
{\it Cocycle Conjugacy Classes of Shifts on the Hyperfinite II$_1$ Factor,
J. Funct. Anal. {\bf 121} (1994), 275-295.

\end